\documentclass[a4paper,DIV=12,9pt]{scrartcl}

\usepackage[utf8]{inputenc}
\usepackage{amsfonts,amssymb,amsmath}
\usepackage{epsfig}
\usepackage{MnSymbol}
\usepackage{booktabs}
\usepackage{cite}
\usepackage{float}
\usepackage{hyperref}

\usepackage{graphics}
\usepackage{subfig}
\usepackage{color}

\usepackage{enumitem}

\newcommand{\ud}{\,\mathrm{d}}
\newcommand{\R}{\mathbb{R}}
\newcommand{\N}{\mathbb{N}}

\renewcommand{\vec}{\boldsymbol}

\numberwithin{equation}{section}
\numberwithin{figure}{section}
\numberwithin{table}{section}

\newtheorem{defi}{Definition}[section]
\newtheorem{thm}[defi]{Theorem}
\newtheorem{lem}[defi]{Lemma}
\newtheorem{rem}[defi]{Remark}

\newtheorem{ass}[defi]{Assumption}
\newtheorem{prob}[defi]{Problem}

\newenvironment{mproof}{\paragraph{Proof.}}{\hfill$\blacksquare$}

\makeatletter
\let\@fnsymbol\@arabic
\makeatother

\begin{document}
	
	\title{Convergence of a continuous Galerkin method for hyperbolic-parabolic systems}
	
	\author{
		Markus Bause$^\ast$\thanks{bause@hsu-hh.de (corresponding author)}\;, 
		Mathias Anselmann$^\ast$, 
		Uwe Köcher$^\ast$,  
		Florin A.\ Radu$^\dag$\\
		{\small ${}^\ast$ Helmut Schmidt University, Faculty of
			Mechanical and Civil Engineering, Holstenhofweg 85,}\\ 
		{\small 22043 Hamburg, Germany}\\
		{\small $^\dag$ University of Bergen,  Center for Modeling of Coupled Subsurface Dynamics,}\\
		{\small Department of Mathematics, All\'egaten 41, 50520  Bergen, Norway}
	}
	
	\date{}
	
	\maketitle

\begin{abstract}
We study the numerical approximation by space-time finite element methods of a multi-physics system coupling hyperbolic elastodynamics with parabolic transport and modeling poro- and thermoelasticity. The equations are rewritten as a first-order system in time. Discretizations by continuous Galerkin methods in time and inf-sup stable pairs of finite element spaces for the spatial variables are investigated. Optimal order error estimates are proved by an analysis in weighted norms that depict the energy of the system's unknowns. A further important ingredient and challenge of the analysis is the control of the couplings terms. The techniques developed here can be generalized to other families of Galerkin space discretizations and advanced models. The error estimates are confirmed by numerical experiments, also for higher order piecewise polynomials in time and space. The latter lead to algebraic systems with complex block structure and put a facet of challenge on the design of iterative solvers. An efficient solution technique is referenced. 
\end{abstract}



\section{Introduction}
\label{Sec:Introduction}

In this work we study the numerical approximation of the coupled equations 
\begin{subequations}
\label{Eq:HPS}
\begin{alignat}{3}
\label{Eq:HPS_1}
\rho \partial_t^2 \vec u - \vec \nabla \cdot (\vec C \vec \varepsilon (\vec u)) + \alpha \vec \nabla p & = \rho \vec f\,, && \quad \text{in } \;
\Omega  \times (0,T]\,,\\[1ex]
\label{Eq:HPS_2}
c_0\partial_t p + \alpha \vec \nabla \cdot \partial_t \vec u  - \vec \nabla \cdot (\vec K \vec \nabla p)  & = g\,, && \quad \text{in } \; \Omega \times (0,T]\,,\\[1ex]
\label{Eq:HPS_3}
 \vec u (0) = \vec u_0\,, \quad \partial_t \vec u (0) = \vec u_1\,, \quad p(0) & = p_0\,, && \quad \text{in } \; \Omega \,,\\[1ex]
\label{Eq:HPS_4}
 \vec u = \vec 0\,, \quad p& =0\,,  && \quad \text{on } \; \partial\Omega \times (0,T]\,.
\end{alignat}
\end{subequations}
Under the below made assumptions about the coefficients of \eqref{Eq:HPS}, this is a system of mixed hyperbolic-parabolic type. It is considered in the open Lipschitz bounded domain $\Omega \subset \R^d$, with $d\in \{2,3\}$, and the time interval $[0,T]$ with some final time $T>0$. For simplicity, Dirichlet boundary conditions are prescribed here in \eqref{Eq:HPS_4}. Important applications of the model \eqref{Eq:HPS}, that is studied as a prototype system, arise in poro- and thermoelasticity. In poroelasticity (cf.~\cite{S00} and \cite{B41,B55,B72}), where Eqs.~\eqref{Eq:HPS} are referred to as the dynamic Biot model, the system \eqref{Eq:HPS} is used to describe flow of a slightly compressible viscous fluid through a deformable porous matrix. The small deformations of the matrix are described by the Navier equations of linear elasticity, and the diffusive fluid flow is described by Duhamel’s equation. The unknowns are the effective solid phase displacement $\vec u$ and the effective fluid pressure $p$. The quantity $\vec \varepsilon (\vec u):= (\vec \nabla u + (\vec \nabla u)^\top)/2$ denotes the symmetrized gradient or strain tensor. Further, $\rho$ is the effective mass density, $\vec C$ is Gassmann’s fourth order effective elasticity tensor, $\alpha$ is Biot’s pressure-storage coupling tensor, $c_0$ is the specific storage coefficient and $\vec K$ is the permeability field. For simplicity, the positive quantities $\rho>0$, $\alpha>0$ and $c_0 >0$ are assumed to be constant in space and time. Moreover, the tensors $\vec C$ and $\vec K$ are assumed to be symmetric and positive definite and independent of the space and time variables as well. In thermoelasticity (cf.~\cite{JR18} and \cite{C72,L86}), the system \eqref{Eq:HPS} describes the flow of heat through an elastic structure. In that context, $p$ denotes the temperature, $c_0$ is the specific heat of the medium, and $\vec K$ is the conductivity. Then, the quantity $\alpha \vec \nabla p$ arises from the thermal stress in the structure, and the term $\alpha \vec \nabla \cdot \partial_t \vec u$  corresponds to the internal heating due to the dilation rate. For the sake of physical realism, the often used uncoupling assumption in which this term is deleted from the diffusion equation is not made here. Well-posedness of \eqref{Eq:HPS} is ensured. For this, we refer to \cite{JR18,S89,STW22} where well-posedness of \eqref{Eq:HPS} is shown by different mathematical techniques, by semigroup methods \cite[Thm.\ 2.2]{JR18}, Rothe's method \cite[Thm.\ 4.18 and Cor.\ 4.33]{S89} and Picard's theorem \cite[Thm.\ 6.2.1]{STW22}. To enhance physical realism, generalizations of the system \eqref{Eq:HPS} are presented in, e.g., \cite{BKNR22,JR18,MW12} and the references therein. 

The coupled hyperbolic-parabolic structure of the system \eqref{Eq:HPS} of partial differential equations adds an additional facet of complexity onto its numerical simulation.  A natural and promising approach for the numerical approximation of coupled systems is given by the application of space-time finite element methods that are based on variational formulations in space and time. Therein, the discrete unknown functions are defined on the entire space-time domain $\Omega\times I$ and can be expanded in terms of finite element basis functions. This facilitates the discretization of even complex coupling terms, for instance, of combined spatial and temporal derivatives or convolution integrals (cf.\ \cite{MW12}). The derivatives in the second of the terms in \eqref{Eq:HPS_2} can be computed naturally, without any further approximation. In this work we are proposing a space-time finite element approximation of the system \eqref{Eq:HPS} by continuous in space and time finite element methods.  For this, the hyperbolic subproblem \eqref{Eq:HPS_1} is rewritten as a first-order system in time. In particular, continuous Galerkin methods provide energy conservative discretizations of wave equations (cf.~\cite[Sec.~6]{BKRS20}), where the energy is measured by $E(t) :=( \| \nabla u(t) \|_{L^2(\Omega)}^2 + \| \partial_t u(t) \|_{L^2(\Omega)}^2 )^{1/2}$ in the scalar-valued case. Thus, continuous Galerkin methods preserve a key structure of solutions to the continuous problem on the discrete level. Here, the continuous Galerkin discretization is considered as a prototype scheme for miscellaneous families of space-time finite element methods. We refer to \cite{AB20,ABBM20,BM21} for the construction of $C^k$-conforming variational time discretizations, for some $k\geq 1$. In this work, we aim to elaborate the treatment of the coupling in \eqref{Eq:HPS} in the error analysis with the perspective of getting optimal order error estimates. We like to present our key arguments and not to overburden the error analysis with the additional terms arisng in discontinuous space discretizations. The error analysis offers the potential and flexibility for its extension to spatial approximations by enriched Galerkin methods (cf.\ \cite{SL09,LLW16,VSBW18}) or discontinuous Galerkin approaches (cf., e.g., \cite{ADMQ16,DF15,GSS06,GS09,K15,PE12}).  Also, for the application of discontinuous Galerkin space discretizations to the quasi-static Biot system, that differs from \eqref{Eq:HPS} by neglecting the acceleration term $\rho \partial_t^2 \vec u$ in \eqref{Eq:HPS_1},  we refer to \cite{PW08,B19}. Unsteady spatial approximations yield appreciable advantages, for instance, for the construction of iterative solver (cf., e.g., \cite{K15}) or the computation by post-processing of locally mass conservative (fluid) fluxes (cf.~\cite{LLW16}) from the variable $p$ of \eqref{Eq:HPS}. The latter is of importance if the system \eqref{Eq:HPS} is coupled further with the transport of species dissolved in the fluid. Discontinuous Galerkin time discretizations (cf.\ \cite{J93,T06}) are not considered here due to their lack of energy conservation for second-order hyperbolic problems. Further, continuous Galerkin methods in time are superior over discontinuous ones regarding the ratio of accuracy, quantified by the convergence rate, over the number of (temporal) degrees of freedom that have to be computed effectively. By an appropriate choice of the trial basis, one temporal degree of freedom is obtained directly by an algebraic relation, which can be exploited to reduce the algebraic system's size; cf., e.g., \cite{HST13}.         
 
 
 The coupling of \eqref{Eq:HPS_1} and  \eqref{Eq:HPS_2} encounters new challenges for the error analysis of numerical schemes and shows a strong link to the mixed approximation by inf-sup stable pairs of finite elements of the Navier--Stokes system; cf.\ \cite{J16}. For this, we note that  \eqref{Eq:HPS} yields a Stokes-type structure for the tuple $(\partial_t \vec u,p)$ in the limit of vanishing coefficients $c_0$ and $\vec K$ such that the well-known stability issues of mixed Stokes approximations emerge and argue either for inf-sup stable pairs of finite element spaces for $\vec u$ and $p$ or for the stabilization of equal-order spatial discretizations. Here, we apply the first of the alternatives and use inf-sup stable pairs of finite element space for the spatial discretization. 
 
 For the approximation of the equations \eqref{Eq:HPS}, rewritten as a first-order system in time with the additional variable $\vec v = \partial_t \vec u$,  by continuous finite element methods of piecewise polynomials of order $k \geq 1$  in time and of order $r\geq 1$ for $p$ as well as of order $r+1$ for $\vec u$ and $\vec v$ in space we show in Thm.~\ref{Thm:MainRes} that the discrete functions $\vec u_{\tau,h} $, $\vec v_{\tau,h} $ and $p_{\tau,h} $ satisfy 
\begin{equation}
\label{Intro:MainEst}
\max_{t\in [0,T]} \big\{\|\nabla (\vec u(t) - \vec u_{\tau,h} (t))\|  + \| \vec v(t)- \vec v_{\tau,h} (t)\|  + \|p(t) - p_{\tau,h}(t)\|\big\} \leq c (\tau^{k+1}+ h^{r+1})\,.
\end{equation}
The error estimate \eqref{Intro:MainEst} is based on energy-type arguments where the energy is measured in a weighted norm. This is essential for the application of the discrete Gronwall inequality. Further, a careful treatment of the coupling terms in \eqref{Eq:HPS} is required to bound their contributions properly which is done here by the choice of suitable test functions along with the application of integration by parts for the time variable. The energy analysis bears out the quantities on left-hand side of \eqref{Intro:MainEst} as its natural errors. Thus, a control of the error in the elastic energy quantity $E(t) = (\| \nabla \vec u(t)\|^2 + \| \partial_t \vec u(t) \|^2)^{1/2}$ of the second-order hyperbolic equation and of the error in the magnitude $\|p(t)\|$ of the unknown of the parabolic subproblem is obtained. Estimate \eqref{Intro:MainEst} is of optimal order with respect to the error quantity $E(t)$ and the pressure $p$. A separation of the errors $\|\nabla (\vec u-\vec u_{\tau,h})\|$ and $\|\vec v-\vec v_{\tau,h}\|$ in their estimation, offering the possibility to increase the spatial convergence order of $\|\vec v-\vec v_{\tau,h}\|$ to $r+2$, does not become feasible by our energy-type arguments. This is due to the fact that \eqref{Eq:HPS_1} is rewritten as a first-order system in time. The error analysis for the resulting system needs test functions that are adapted to this mixed structure of partial and ordinary differential equations; cf.\  Rem.~\ref{Rem:OptOrder}. Thereby, decoupling mechanisms are inhibited.

The continuous in time Galerkin discretization is known to be superconvergent in the temporal nodes, more precisely, in the Gauss--Lobatto quadrature points of the subintervals of the time mesh, if $k\geq 2$. For the heat and wave equation, superconvergence is studied in \cite{AM89} and \cite{BKRS20}, respectively, and for systems of ordinary differential equations in \cite{BM21}. We conjecture and show numerically that 
\begin{equation}
	\label{Intro:MainEst2}
	\max_{n=1,\ldots , N} \big\{\|\nabla (\vec u(t_n) - \vec u_{\tau,h} (t_n))\|  + \| \vec v(t_n)- \vec v_{\tau,h} (t_n)\|  + \|p(t) - p_{\tau,h}(t)\|\big\} \leq c (\tau^{2k}+ h^{r+1})\,,
\end{equation}
is satisfied. A proof of \eqref{Intro:MainEst2} remains an open problem and is left as a work for the future. Here, we firstly prove \eqref{Intro:MainEst} that is expected to be a prerequisite for showing the result of superconvergence \eqref{Intro:MainEst2}. 

This work is organized as follows. In Sec.~\ref{Sec:Not}, notations and auxiliary results are introduced. In Sec.~\ref{Sec:Scheme}, our approximation of \eqref{Eq:HPS} is presented. In Sec.~\ref{Sec:ErrAna}, the error estimation is done and \eqref{Intro:MainEst} is proved. Finally, in Sec.~\ref{Sec:NumExp}, the results of our numerical experiments are summarized. An efficient iterative solver for the arising algebraic system is referenced. 

\section{Notations, finite element spaces and auxiliaries}
\label{Sec:Not}

\subsection{Notations}

In this work, standard notation is used. We denote by $H^m(\Omega)$ the Sobolev space of $L^2(\Omega)$ functions with weak derivatives up to order $m$ in $L^2(\Omega)$. We let $H^1_0(\Omega)=\{u\in H^1(\Omega) \mid u=0 \mbox{ on } \partial \Omega\}$. For short, we skip the domain $\Omega$ in the notation. Thus, we put $L^2=L^2(\Omega)$, $H^m=H^m(\Omega)$ and $H^1_0=H^1_0(\Omega)$. By $H^{-1}=H^{-1}(\Omega)$ we denote the dual space of $H^1_0$. For vector-valued functions we write those spaces bold. By $\llangle \cdot, \cdot \rrangle$ we define the $L^2$ inner product on the product space $(L^2)^2$. For the norms of the Sobolev spaces the notation is 
\begin{align*}
	\| \cdot \| := \| \cdot\|_{L^2}\,,\qquad 
	\| \cdot \|_m := \| \cdot\|_{H^m}, \,\, \mbox{ for } m \in \N_0\,, \;\; (H^0:=L^2)\,.
\end{align*} 
For a Banach space $B$ we let $L^2(0,T;B)$, $C([0,T];B)$ and $C^m([0,T];B)$, $m\in\N$, be the Bochner spaces of 
$B$-valued functions, equiped with their natural norms. Further, for a subinterval $J\subseteq [0,T]$, we will use the notations $L^2(J;B)$, $C^m(J;B)$ and $C^0(J;B):= C(J;B)$ for the corresponding Bochner spaces.

In what follows, the constant $c$ is generic and indepedent of the size of the space and time meshes. The value of $c$ can depend on 
norms of the solution to \eqref{Eq:HPS}, the regularity of the space mesh, the polynomial degrees used for the space-time discretization and the data (including $\Omega$).

\subsection{Finite element spaces}
\label{Subsec:FES}

For the time discretization, we decompose the time interval $I=(0,T]$ into $N$ subintervals $I_n=(t_{n-1},t_n]$, $n=1,\ldots,N$, where $0=t_0<t_1< \cdots < t_{N-1} < t_N = T$ such that $I=\bigcup_{n=1}^N I_n$. We put $\tau := \max_{n=1,\ldots, N} \tau_n$ with $\tau_n = t_n-t_{n-1}$. Further, the set $\mathcal{M}_\tau := \{I_1,\ldots,
I_N\}$ of time intervals is called the time mesh. For a Banach space $B$ and any $k\in
\N_0$, we let 
\begin{equation}
\label{Def:Pk}
	\mathbb P_k(I_n;B) := \bigg\{w_\tau \,: \,  I_n \to B \,, \; w_\tau(t) = \sum_{j=0}^k 
	W^j t^j \; \forall t\in I_n\,, \; W^j \in B\; \forall j \bigg\}\,.
\end{equation}
For an integer $k\in \N$, we introduce the space
\begin{equation}
	\label{Eq:DefXk} 
	X_\tau^k (B) := \left\{w_\tau \in C(\overline{I};B) \mid w_\tau{}_{|I_n} \in
	\mathbb P_k(I_n;B)\; \forall I_n\in \mathcal{M}_\tau \right\}
\end{equation}
of globally continuous in time functions and for an integer $l\in \N_0$ the space
\begin{equation}
	\label{Eq:DefYk}
	Y_\tau^{l} (B) := \left\{w_\tau \in L^2(I;B) \mid w_\tau{}_{|I_n} \in
	\mathbb P_{l}(I_n;B)\; \forall I_n\in \mathcal{M}_\tau \right\}
\end{equation}
of global $L^2$-functions in time. For a function $w:I\to B$ that is piecewise continuous with respect to the time mesh $\mathcal{M}_{\tau}$, we define by
\begin{equation}
\label{Eq:Defw_In_bdr}
w (t_n^+):= \lim_{t\to t_n+0} w(t)
\qquad\text{and}\qquad
w(t_n) := \lim_{t\to t_n-0} w(t)
\end{equation}
the one-sided limits of $w$. For brevity, we skip the upper index for the argument of $w$ in the second of the definitions, since by definition $I_n = (t_{n-1},t_n]$ such that $w_{|I_n}(t_n)$ is well-defined.  

For the space discretization, let $\mathcal{T}_h=\{K\}$ be a family of shape-regular meshes of $\Omega$ consisting of quadrilateral or hexahedral elements $K$ with mesh size $h>0$ that we use for our computations (cf.\ Sec.~\ref{Sec:NumExp}). Further, for any $r\in \N$ let $V_h^r$ be the finite element space that is built on the mesh of quadrilateral or hexahedral elements and is given by 
\begin{equation}
\label{Eq:DefVh}
V_h^r := \left\{v_h \in C(\overline \Omega) \mid v_h{}_{|K}\in \mathbb Q_r(K) \, \forall K  \in \mathcal T_h \right\}\cap H^1_0(\Omega)\,,
\end{equation}
where $\mathbb Q_r(K)$ is the space defined by the reference mapping of polynomials on the reference element with maximum degree $r$ in each variable. For vector-valued functions we write the space $V_h^r$ bold.

\subsection{Auxiliaries: Quadrature formulas and interpolation operators in time}

We will need some quadrature formulas and interpolation and projection operators acting on the time variable. For the continuous in time finite element method, a natural choice is to consider the $(k+1)$-point Gau{ss}--Lobatto quadrature formula on each time interval $I_n=(t_{n-1},t_n]$,
\begin{equation}
	\label{Eq:GLF}
	Q_n(w) := \frac{\tau_n}{2}\sum_{\mu=0}^{k} \hat \omega_\mu w{}_{|I_n}(t_{n,\mu}) \approx 
	\int_{I_n} w(t) \ud t\,,  
\end{equation}
where $t_{n,\mu}=T_n(\hat t_{\mu})$, for $ \mu = 0,\ldots ,k$, are the quadrature points on $\bar I_n$ and $\hat  \omega_\mu$ the corresponding weights. Here, $T_n(\hat t):=(t_{n-1}+t_n)/2 + (\tau_n/2)\hat t$ is the affine transformation from the reference interval $\hat I = [-1,1]$ to $I_n$ and $\hat t_{\mu}$, for $\mu = 0,\ldots,k$, are the Gau{ss}--Lobatto quadrature points on $\hat I$. We note that for the Gau{ss}--Lobatto formula the identities $t_{n,0}=t_{n-1}$ and $t_{n,k}=t_{n}$ are satisfied and that the values $w_{|I_n}(t_{n,\mu})$ for $\mu\in\{0,k\}$ denote the corresponding one-sided limits of values $w(t)$ from the interior of $I_n$ (cf.\ \eqref{Eq:Defw_In_bdr}). It is known that formula \eqref{Eq:GLF} is 
exact for all polynomials in $\mathbb P_{2k-1} (I_n;\R)$. For the Gau{ss}--Lobatto quadrature points $t_{n,\mu}$, with $n=1,\ldots, N$ and $\mu = 0,\ldots ,k$, we also define the global Lagrange interpolation operator $I_\tau :C^0(\overline I;L^2)\mapsto X^{k}_\tau(L^2)$ by means of
\begin{equation}
	\label{Eq:DefLagIntOp}
	I_\tau w(t_{n,\mu}) = w(t_{n,\mu})\,, \quad \mu=0,\ldots,k\,, \; n=1,\ldots, N\,.
\end{equation}

The $k$-point Gau{ss} quadrature formula on $I_n$ is denoted by 
\begin{equation}
	\label{Eq:GF}
	Q_n^{\operatorname G}(w) := \frac{\tau_n}{2}\sum_{\mu=1}^{k} \hat 
	\omega_\mu^{\operatorname 
		G} w(t_{n,\mu}^{\operatorname G}) \approx \int_{I_n} w(t) \ud t \,,
\end{equation}
where $t_{n,\mu}^{\text G}=T_n(\hat t_{\mu}^{\,\operatorname G})$, for $\mu =1,\ldots, k$, are the Gauss quadrature 
points on $I_n$ and $\hat \omega_\mu^{\operatorname G}$ the corresponding weights, with $\hat t_{\mu}^{\,\operatorname G}$, for $\mu = 1,\ldots,k$, being the Gau{ss} quadrature 
points on $\hat I$. Formula \eqref{Eq:GF} is also exact for all polynomials in $\mathbb P_{2k-1} (I_n;\R)$. For $n=1,\ldots, N$, the local interpolant $I_{\tau,n}^{\text G} :C^0(\overline I_n;L^2)\mapsto \mathbb P_{k-1}(\overline I_n;L^2)$ is defined by means of 
\begin{equation}
	\label{Eq:DefGaussIntOp}
	I_{\tau,n}^{\text G} w(t_{n,\mu}^{\text G}) = w(t_{n,\mu}^{\text G})\,. \quad \mu=1,\ldots,k\,.
\end{equation}

Further, for a given function $w\in L^2(I;B)$, we define the interpolate
$\Pi^{k-1}_\tau w\in Y^{k-1}_\tau(B)$ such that its restriction $\Pi^{k-1}_\tau w{}_{|I_n}\in\mathbb P_{k-1}(I_n;B)$, $n=1,\ldots, N$, is determined 
by local  $L^2$-projection in time, i.e.  
\begin{equation}
	\label{Def:Pi}
	\int_{I_n} \langle \Pi^{k-1}_\tau  w , q \rangle \ud t  =
	\int_{I_n} \langle  w , q \rangle \ud t 
	\qquad\forall\, q\in \mathbb P_{k-1}(I_n;B) \,.
\end{equation}

\begin{rem}
All operators, that act on the temporal variable only, are applied componentwise to a vector field $\vec F= (F_0,\ldots,F^d)^\top$, for instance, $I_\tau \vec F = (I_\tau F_0, \ldots, I_\tau F_d)^\top$. This is tacitly used below. 
\end{rem}

The following result (cf.\ \cite[Eq.~(2.6)]{KM04} and \cite[Lem.\ 4.5]{BKRS20}) is proved easily.

\begin{lem}
\label{Lem:PropGF} 
Consider the Gau{ss} quadrature formula \eqref{Eq:GF}. For all $n=1,\ldots,N$ there holds that 
\begin{subequations}
\label{Eq:PropGF} 
\begin{alignat}{2}	
	\label{Eq:PropGF_0} 
	\Pi_\tau^{k-1} w(t) & = I_{\tau,n}^{\text G} w(t) \,, & & \quad \text{for} \;\; t\in I_n\,,\\[1ex]
\label{Eq:PropGF_1} 
\Pi^{k-1}_\tau  w(t_{n,\mu}^{\operatorname G}) & = w(t_{n,\mu}^{\operatorname G})\,, &&  \quad \text{for} \;\; \mu = 
		1,\ldots ,k\,,
\end{alignat}
\end{subequations}
for all polynomials $w\in \mathbb P_{k}(I_n;L^2)$.
\end{lem}

Finally, we recall the following $L^\infty$--$L^2$ inverse inequality; cf.~\cite[Eq.\ (2.5)]{KM04}.
\begin{lem}
\label{Lem:LinfL2} 
For all $n=1,\ldots,N$ there holds that
\begin{equation}
\label{Eq:LinfL2}
\| w \|_{L^{\infty}(I_n;\R)} \leq c \tau_n^{-1/2} \| w \|_{L^2(I_n;\R)} 
\end{equation} 
for all polynomials $w\in \mathbb P_{k}(I_n;\R)$.
\end{lem}

\section{The fully discrete scheme and preparation for the error analysis}
\label{Sec:Scheme}

Here we propose our discretization of \eqref{Eq:HPS} by continuous finite element methods in time and space. For the discretization we rewrite Eq.~\eqref{Eq:HPS_1} as a first-order system in time such that time-discretization schemes designed for first-order systems of ordinary differential equations become applicable.  For this, we put $\vec v := \partial_t \vec u$. 

\subsection{Bilinear forms and discrete operators}

Here we introduce (bi-)linear forms for the discrete variational formulation and further operators related to the spatial discretization. For $\vec u, \vec v, \vec \phi\in \vec H^1_0$, $p, \psi \in H^1_0$, $\vec f \in \vec H^{-1}$ and $g\in H^{-1}$ we put 
\begin{align*}
A(\vec u, \vec \phi) & :=  \langle \vec C \vec \epsilon(\vec u), \vec \varepsilon(\vec \phi) \rangle\,, &
B(p, \psi) & := \langle \vec K \nabla p, \nabla \psi \rangle\,, & C(\vec v,\psi)  & :=  - \alpha \langle \nabla \cdot \vec v, \psi \rangle\,, \\[1ex] F(\vec \phi) & := \langle \rho \vec f, \vec \phi\rangle \,, & G(\psi) & := \langle g, \psi \rangle\,.
\end{align*}
 
Firstly, we address the discretization of the hyperbolic equation \eqref{Eq:HPS_1}. By $\vec P_h: \vec L^2\mapsto \vec V_h^{r+1}$ we denote the $\vec L^2$-orthogonal projection onto $\vec V_h^{r+1}$ such that, for 
$\vec w\in \vec L^2$, the identity 
\begin{equation*}
	\langle \vec P_h \vec w, \vec \phi_h \rangle = \langle \vec w, \vec \phi_h\rangle  
\end{equation*}
is satisfied for all $\vec \phi_h\in \vec V_h^{r+1}$. The operator $\vec R_h: \vec H^1_0 \mapsto \vec V_h^{r+1}$ defines the elliptic 
projection onto $\vec V_h^{r+1}$ such that 
\begin{equation}
	\label{Def:vRh}
	\langle  \vec C \vec \varepsilon (\vec R_h \vec w), \vec \varepsilon(\vec \phi_h) \rangle = \langle  \vec C \vec \varepsilon (\vec w), \vec \varepsilon(\vec \phi_h) \rangle
\end{equation}
for $\vec w\in \vec H^1_0$ and all $\vec \phi_h\in \vec V_h^{r+1}$. We let $\vec A_h: \vec H^1_0 \mapsto \vec V_h^{r+1}$ be the discrete operator that is defined by 
\begin{equation}
\label{Eq:DefAh}
\langle \vec A_h \vec w , \vec \phi_h \rangle = A(\vec w,\vec \phi_h)
\end{equation}
for all $\vec \phi_h\in \vec V_h^{r+1}$. Then, for $\vec w \in \vec H^1_0 \cap \vec H^2$ it holds that 
\begin{equation}
	\label{Eq:AhwAw}
	\langle \vec A_h \vec w , \vec \phi_h \rangle = \langle \vec C \vec \epsilon(\vec w), \vec \varepsilon(\vec \phi_h) \rangle\ = \langle \vec A \vec w , \vec \phi_h \rangle 
\end{equation}
for $\vec \phi_h\in \vec V_h^{r+1}$, where $\vec A: \vec H^1_0 \rightarrow \vec H^{-1}$ is defined by $\langle \vec A\vec w, \vec \phi\rangle := A(\vec w, \vec \phi)$ for $\vec \phi \in \vec H^1_0$. Thus, $\vec A_h \vec w = \vec P_h \vec A \vec w $ for $\vec w\in \vec H^1_0 \cap \vec H^2$.

Further, let $\vec{\mathcal L}_h:  \vec H^1_0 \times \vec L^2 \mapsto \vec V_h^{r+1} \times \vec V_h^{r+1}$ be 
defined by 
\begin{equation}
	\label{Eq:DefLh} 
	\vec{\mathcal L}_h := \begin{pmatrix}
		\vec 0 & -\vec I\\ \vec A_h & \vec 0 
	\end{pmatrix}\,.
\end{equation}
Then, for $\vec U= (\vec U_1,\vec U_2)\in (\vec H^1_0\cap \vec H^2)\times \vec L^2$ we have that 
\begin{equation*}
	\llangle \vec{\mathcal L}_h \vec U, \vec \Phi_h \rrangle = \langle -\vec U_2 , \vec \Phi_h^1\rangle + \langle \vec C
	\vec \varepsilon(\vec U_1), \vec \varepsilon (\vec \Phi_h^2) \rangle = \langle - \vec U_2, \vec \Phi_h^1\rangle + \langle \vec A 
	\vec U_1, \vec \Phi_h^2\rangle = \llangle \vec{\mathcal L} \vec U, \vec \Phi_h \rrangle
\end{equation*}
for $\vec \Phi_h = (\vec \Phi_h^1,\vec \Phi_h^2)^\top \in \vec V_h^{r+1}\times \vec V_h^{r+1}$, where $\vec{\mathcal L}: \vec H^1_0 \times \vec L^2 \rightarrow \vec L^2\times \vec H^{-1} $, with $D(\vec{\mathcal L})=\vec H^1_0 \times \vec L^2$, is defined by
$\vec{\mathcal L} := \begin{pmatrix}
	\vec 0 & -\vec I\\ \vec A & \vec 0 
\end{pmatrix}. $

Secondly, we address the discretization of the parabolic equation~\eqref{Eq:HPS_2}. By $P_h: L^2\mapsto V_h^r$ we denote the $L^2$-orthogonal projection onto $V_h^r$ such that, for 
$w\in L^2$, the identity 
\begin{equation*}
	\langle P_h w, \psi_h \rangle = \langle w, \psi_h\rangle  
\end{equation*}
is satisfied for all $\psi_h\in V_h^r$. The operator $R_h: H^1_0 \mapsto V_h^r$ defines the elliptic 
projection onto $V_h^r$ such that, for $w\in H^1_0$,
\begin{equation}
	\label{Def:Rh}
	\langle  \vec K \nabla R_h w, \nabla \psi_h \rangle = \langle  \vec K \nabla w, \nabla \psi_h \rangle   
\end{equation}
for all $\psi_h\in V_h^r$. Let $B_h: H^1_0  \mapsto V_h^r$ be the discrete operator that is defined by 
\begin{equation}
	\label{Eq:DefBh}
	\langle B_h w , \psi_h \rangle :=  B(w, \psi_h)   
\end{equation}
for all $v_h\in V_h$. Then, for $w \in H^1_0 \cap H^2$ it holds that 
\begin{equation*}
	\langle B_h w , v_h \rangle = \langle \vec K \nabla w, \nabla v_h \rangle = \langle Bw 
	,v_h\rangle 
\end{equation*}
for all $v_h\in V_h^r$, where $B: H^1_0 \rightarrow H^{-1}$ is defined by $\langle Bw,\psi\rangle  = B(w,\psi)$ for $\psi\in H^1_0$. Thus, $B_h w = P_h Bw $ for $w\in H^1_0 \cap H^2$.

\begin{rem}
We note that discrete functions of $\vec V_{h}^{r+1}$, with some $r\in \N$, will be used for the approximation of the vectorial variable $\vec u$ and discrete functions of $V_{h}^{r}$ for the approximation of the scalar variable $p$; cf.\ Subsec.~\ref{Subsec:CGD}. The projection and discrete differential operators are thus defined for finite element spaces of different polynomial degrees, which is not expressed explicitly by the notation for brevity. 
\end{rem}

\subsection{Continuous Galerkin discretization}
\label{Subsec:CGD}

Here, we formulate our space-time finite element approximation of the system \eqref{Eq:HPS}. For the discretization in time, the continuous Galerkin method is applied; cf.\ \cite{AM89,BKRS20,BRK17,FP96,KM04,S10}. Precisely, the time discretization is of Petrov--Galerkin type. For the discretization in space, a continuous finite element approach, based on inf-sup stable pairs of finite elements, is used. Generalization of the error analysis to other families of Galerkin space discretizations that offer appreciable advantages, like local mass conservation, appear feasible. Restricting ourselves to the family of Taylor--Hood pairs of finite element spaces is done in order to carve out the key arguments of our error analysis. 


We make the following assumption about the discrete initial values $\vec u_{0,h},\vec v_{0,h}\in \vec V_h^{r+1}$ and $p_{0,h}\in V_h^r$.

\begin{ass}
\label{AssIV}
Let $\vec u_{0,h},\vec v_{0,h}\in \vec V_h^{r+1}$ and $p_{0,h}\in V_h^r$ be chosen such the approximation properties 
\begin{subequations}
\begin{alignat}{2}
\label{Eq:AssIV0}
\| \nabla (\vec R_h \vec u_0 - \vec u_{0,h})\| & \leq c h^{r+1}\| \vec u_0 \|_{r+2}\,,\\[1ex]
\label{Eq:AssIV1}
\| \vec R_h \vec u_1 - \vec v_{0,h}\| & \leq c h^{r+2} \| \vec u_1 \|_{r+2}\,,\\[1ex]
\label{Eq:AssIV2}
\| R_h p_0 - p_{0,h}\| & \leq c h^{r+1}\|p_0\|_{r+1}
\end{alignat}
\end{subequations}
are satisfied for $\vec u_0,\vec u_1 \in \vec H^1_0\cap \vec H^{r+2}$ and $p_0\in H^1_0\cap H^{r+1}$, where $\vec R_h$ and $R_h$ are defined by \eqref{Def:vRh} and \eqref{Def:Rh}, respectively. 
\end{ass}

We use a temporal test basis that is supported on the subintervals $I_n$; cf.\ \cite{S10,BKRS20}. Then, a time marching process is obtained. In that, we assume that the trajectories $\vec u_{ \tau,h}$, $\vec v_{ \tau,h}$ and $p_{ \tau,h}$ have been computed before for all $t\in [0,t_{n-1}]$, starting with approximations  $\vec u_{\tau,h}(t_0) :=\vec u_{0,h}$, $\vec v_{\tau,h}(t_0) :=\vec v_{0,h}$ and $p_{\tau,h}(t_0) := p_{0,h}$ of the initial values $\vec u_0$, $\vec u_1$ and $p_{0}$. Then, we consider solving the following local problem on $I_n$.

\begin{prob}[Variational form of $I_n$ problem]
\label{Prob:DPV}
Let $k,r\geq 1$. For given $\vec u_{\tau,h}^{n-1}:=\vec u_{\tau,h}(t_{n-1})\in \vec V_h^{r+1}$, $\vec v_{\tau,h}^{n-1}:=\vec v_{\tau,h}(t_{n-1})\in \vec V_h^{r+1}$, $p_{\tau,h}^{n-1}:= p_{\tau,h}(t_{n-1})\in V_h^{r}$ with $\vec u_{\tau,h}(t_0) :=\vec u_{0,h}$, $\vec v_{\tau,h}(t_0) :=\vec v_{0,h}$ and $p_{\tau,h}(t_0) := p_{0,h}$, find $\vec U_{\tau,h}= (\vec u_{\tau,h},\vec v_{\tau,h})^\top \in (\mathbb P_k(I_n;V_{h}^{r+1}))^d\times (\mathbb P_k(I_n;V_{h}^{r+1}))^d$ and ${p_{\tau,h} \in \mathbb P_k(I_n;V_{h}^{r})}$ such that $\vec U_{\tau,h}(t_{n-1})=(\vec u_{\tau,h}^{n-1},\vec v_{\tau,h}^{n-1})^\top$, $p_{\tau,h}(t_{n-1}) = p_{\tau,h}^{n-1}$ and 
\begin{subequations}
	\label{Eq:DPV_0}
	\begin{align}
		\label{Eq:DPV_1}
		\int_{I_n} \llangle \vec D \partial_t \vec U_{\tau,h} , \vec \Phi_{\tau,h} \rrangle + \llangle \vec{\mathcal L}_h \vec U_{\tau,h}, \vec \Phi_{\tau,h} \rrangle  - \alpha \langle p_{\tau,h},\nabla \cdot \vec \Phi^2_{\tau,h}\rangle \ud t & = Q_n \Big(\llangle \vec F, \vec \Phi_{\tau,h}\rrangle\Big)\,,\\[1ex]
		\label{Eq:DPV_2}
		\int_{I_n} \langle c_0 \partial_t p_{\tau,h},\psi_{\tau,h} \rangle + \alpha \langle \nabla \cdot \partial_t \vec u_{\tau,h},\psi_{\tau,h}\rangle + \langle B_h p_{\tau,h}, \psi_{\tau,h}\rangle \ud t  & = Q_n \Big( \langle g,\psi_{\tau,h}\rangle \Big)
	\end{align}
\end{subequations}
for all $ \vec \Phi_{\tau,h}=(\vec \Phi^1_{\tau,h},\vec \Phi^2_{\tau,h})^\top\in (\mathbb P_{k-1}(I_n;V_{h}^{r+1}))^d\times (\mathbb P_{k-1}(I_n;V_{h}^{r+1}))^d$ and ${\psi_{\tau,h} \in \mathbb P_{k-1}(I_n;V_{h}^{r})}$, where $\vec D := \begin{pmatrix} \vec I_d & \vec 0 \\ \vec 0 & \rho \vec I_d \end{pmatrix}$ with the identity matrix $\vec I_d\in \R^{d,d}$ and $\vec F :=(\vec 0^\top,\rho \vec f^\top)^\top$. 

\end{prob}

\begin{rem}
By means of Lem.~\ref{Lem:Relwj}, given below, the term $\int_{I_n} \langle \nabla \cdot \partial_t \vec u_{\tau,h},\psi_{\tau,h}\rangle \ud t$ in \eqref{Eq:DPV_2} can equivalently be replaced by $\int_{I_n} \langle \nabla \cdot \vec v_{\tau,h},\psi_{\tau,h}\rangle \ud t$.
\end{rem}

By the exactness of the Gauss--Lobatto quadrature formula \eqref{Eq:GLF} for all polynomials in $\mathbb P_{2k-1} (I_n;\R)$ we can recover the variational problem \eqref{Eq:DPV_0} in the following numerically integrated form. 

\begin{prob}[Quadrature form of $I_n$ problem]
\label{Prob:DPQ}
Let $k,r\geq 1$. For given $\vec u_{\tau,h}^{n-1}:=\vec u_{\tau,h}(t_{n-1})\in \vec V_h^{r+1}$, $\vec v_{\tau,h}^{n-1}:=\vec v_{\tau,h}(t_{n-1})\in \vec V_h^{r+1}$, $p_{\tau,h}^{n-1}:= p_{\tau,h}(t_{n-1})\in V_h^{r}$ with $\vec u_{\tau,h}(t_0) :=\vec u_{0,h}$, $\vec v_{\tau,h}(t_0) :=\vec v_{0,h}$ and $p_{\tau,h}(t_0) := p_{0,h}$, find $\vec u_{\tau,h}  \in( \mathbb P_k(I_n;V_{h}^{r+1}))^d$, $\vec v_{\tau,h}  \in (\mathbb P_k(I_n;V_{h}^{r+1}))^d$ and ${p_{\tau,h} \in \mathbb P_k(I_n;V_{h}^{r})}$ such that $\vec u_{\tau,h}(t_{n-1})=\vec u_{\tau,h}^{n-1}$, $\vec v_{\tau,h}(t_{n-1}) = \vec v_{\tau,h}^{n-1}$, $p_{\tau,h}(t_{n-1}) = p_{\tau,h}^{n-1}$ and 
\begin{subequations}
\label{Eq:DPQ_0}
\begin{align}
\label{Eq:DPQ_1}
Q_n \big(\langle \partial_t \vec u_{\tau,h} , \vec \phi_{\tau,h} \rangle  - \langle \vec v_{\tau,h} , \vec \phi_{\tau,h} \rangle \big) & = 0\,,\\[1ex]
\label{Eq:DPQ_2}
Q_n \Big(\langle \rho \partial_t \vec v_{\tau,h} , \vec \chi_{\tau,h} \rangle + A(\vec u_{\tau,h}, \vec \chi_{\tau,h} ) + C(\vec \chi_{\tau,h},p_{\tau,h})\Big) & = Q_n \Big(F(\vec \chi_{\tau,h})\Big)\,,\\[1ex]
\label{Eq:DPQ_3}
Q_n \Big(\langle c_0 \partial_t p_{\tau,h},\psi_{\tau,h} \rangle  - C(\partial_t \vec u_{\tau,h},\psi_{\tau,h})+ B(p_{\tau,h}, \psi_{\tau,h})\Big) & = Q_n \Big( G(\psi_{\tau,h})\Big)
\end{align}
\end{subequations}
for ${\vec \phi_{\tau,h} \in (\mathbb P_{k-1}(I_n;V_{h}^{r+1}))^d}$, ${\vec \chi_{\tau,h} \in (\mathbb P_{k-1}(I_n;V_{h}^{r+1}))^d}$ and ${\psi_{\tau,h} \in \mathbb P_{k-1}(I_n;V_{h}^{r})}$.	
\end{prob}

\begin{rem}
\label{Rem:NEOA}
\begin{itemize}
\item Problem~\ref{Prob:DPV} or \ref{Prob:DPQ}, respectively, yields a globally continuous in time discrete solution 
\begin{equation*}
\label{Eq:GlobDiscSol} 
(\vec u_{\tau,h},\vec v_{\tau,h},p_{\tau,h}) \in (X_\tau^{k}(V_{h}^{r+1}))^d \times (X_\tau^{k}(V_{h}^{r+1}))^d \times X_\tau^{k}(V_{h}^{r})\,.
\end{equation*}

\item A non-equal order spatial approximation of the unknowns $(\vec u,p)$ in the spaces $\vec V_{h}^{r+1}\times V_{h}^{r}$, built from the Taylor--Hood pair of element spaces, is applied here. The inf-sup (or LBB) stability condition is satisfied by this choice of spaces; cf.~\cite{J16}. For vanishing coefficients $c_0\rightarrow 0$ and $\vec K \rightarrow \vec 0$, a Stokes-type system structure is obtained in \eqref{Eq:HPS} for the variables $\partial_t \vec u$ and $p$ such that the well-known stability issues of mixed approximations of the Stokes system emerge in the limit case of vanishing $c_0$ and $\vec K$; cf.~\cite{J16}. Therefore, equal order spatial discretizations do not become feasible without any additional stabilization of the discretization. For a more detailed discussion of stability properties for the quasi-static Biot system we also refer to, e.g., \cite{ML92,ML94,MTL96,RHOAGZ17}. 

\item In Problem \ref{Prob:DPQ}, the Gauss-Lobatto quadrature formula is applied. This allows an efficient implementation of the continuity constraints at the discrete time nodes $t_n$, for $n= 0,\ldots, N-1$, in computer codes (cf.\ \cite{KB14,K15}) and, thus, is the most natural approach for the continuous Galerkin approximation in time. In the error analysis, the Gauss quadrature formula \eqref{Eq:GF}, that is also exact for all polynomials in $\mathbb P_{2k-1} (I_n;\R)$, is used as well.

\end{itemize}
\end{rem}

\subsection{Preparation for the error analysis}

Here we present some auxiliaries that will used below in the error analysis. Firstly, we introduce some special approximation $\vec w = (\vec w_1,\vec w_2)$ of the solution $(\vec u, \vec v)$, with $\vec v := \partial_t \vec u$, that has been defined in \cite{KM04}.
\begin{defi}[Special approximation $(\vec w_1,\vec w_2)$ of $(\vec u, \partial_t \vec u)$]
\label{Def:W}
Let $\vec u\in C^1(\overline I;\vec H^1_0)$ be given. On $I_n=(t_{n-1},t_n]$ we define 
\begin{equation}
\label{Eq:DefW}
\vec w_1 := I_\tau \Big( \int_{t_{n-1}}^t \vec w_2(s)\ud s + \vec R_h \vec u(t_{n-1})  \Big)\,,
\qquad\text{where}\qquad  \vec w_2 := I_\tau (\vec R_h\partial_t \vec u) \,.
\end{equation}
Further, we put $\vec w_1(0) := \vec R_h \vec u(0)$. 	
\end{defi}

In Def.~\ref{Def:W} we simply write $\vec w_j$, for $j=1,2$, instead of $\vec w_j{}_|{}_{I_n}$. The Lagrange interpolation operator $I_\tau$ for the Gau{ss}-Lobatto quadrature points (cf.\ \eqref{Eq:DefLagIntOp}) acts locally on $\bar I_n$ as $I_\tau: C^0(\bar I_n;B)\mapsto \mathbb P_k(I_n;B)$ for any Banach space $B$. The approximations $\vec w_j\in \mathbb (P_{k}(I_n;V_h^{r+1}))^d$, for $j=1,2$, satisfy the following variational equation (cf.\ \cite[Lem.~3.1]{KM04}).  

\begin{lem}
\label{Lem:Relwj}
For $\vec w_1$ and $\vec w_2$, defined in Def.~\ref{Def:W}, there holds for all $\vec \phi_{\tau,h} \in (\mathbb P_{k-1}(I_n;V_{h}^{r+1}))^d$ that 
\begin{equation}
\label{Eq:Relwj}
\int_{I_n} \langle \partial_t \vec w_1,  \vec \phi_{\tau,h}\rangle \ud t = \int_{I_n} \langle \vec w_2,  \vec \phi_{\tau,h}\rangle \ud t\,.  	
\end{equation}
\end{lem}

Further, we need the following auxiliary result for the error analysis.  
\begin{lem}
\label{Lem:upev}
For $\vec y_{\tau,h}, \vec z_{\tau,h} \in (\mathbb P_k(I_n;V_{h}^{r+1}))^d$ let 
\begin{equation}
\label{Eq:Aux01}
\int_{I_n} \langle \partial_t \vec y_{\tau,h}, \vec \phi_{\tau,h}  \rangle - \langle \vec z_{\tau,h}, \vec \phi_{\tau,h} \rangle \ud t = 0	
\end{equation}
be satisfied for all  ${\vec \phi_{\tau,h} \in (\mathbb P_{k-1}(I_n;V_{h}^{r+1}))^d}$. Then, there holds that 	
\begin{equation}
\label{Eq:Aux02}
\partial_t \vec y_{\tau,h}(t_{n,\mu}^{\text G})= \vec z_{\tau,h}(t_{n,\mu}^{\text G})
\end{equation}
for $\mu = 1,\ldots , k$, where $\{t_{n,\mu}^{\text G}\}_{\mu = 1}^k$ are the Gauss quadrate nodes (cf.~\eqref{Eq:GF}) of the subinterval $I_n$.
\end{lem}

\begin{mproof}
Let $l\in \{1,\ldots,k\}$ be arbitrary but fixed and $\vec \phi_{\tau,h}\in (\mathbb P_{k-1}(I_n;V_{h}^{r+1}))^d$ be chosen as  
\begin{equation*}
\vec \phi_{\tau,h}(t) := \xi_n(t) \vec \phi_h \quad \text{with}\quad \xi_n(t) := \prod_{i=1 \atop i \neq l}^k (t-t_{n,i}^{\text G})\in \mathbb P_{k-1}(I_n;\R)\,, \;\; \vec \phi_h \in \vec V_h^{r+1}\,, 
\end{equation*}
and the Gauss quadrature nodes $t_{n,\mu}^{\text G}$, for $\mu = 1,\ldots , k$; cf.~\eqref{Eq:GF}. By the exactness of the Gauss quadrature formula \eqref{Eq:GF} for all polynomials in $\mathbb P_{2k-1}(I_n;\R)$ we deduce from \eqref{Eq:Aux01} that
\begin{align*}
0 & = \int_{I_n} \langle \partial_t \vec y_{\tau,h}, \vec \phi_{\tau,h}  \rangle - \langle \vec z_{\tau,h}, \vec \phi_{\tau,h} \rangle \ud t =	\frac{\tau_n}{2}\sum_{\mu=1}^k \hat \omega_\mu^{\text G}(\langle \partial_t \vec y_{\tau,h}(t_{n,\mu}^{\text G}), \vec \phi_{\tau,h}(t_{n,\mu}^{\text G})  \rangle - \langle \vec z_{\tau,h}(t_{n,\mu}^{\text G}), \vec \phi_{\tau,h}(t_{n,\mu}^{\text G}) \rangle)\\[1ex]
& = \frac{\tau_n}{2} \hat \omega_\mu^{\text G} \xi_{n}(t_{n,l}^{\text G}) (\langle \partial_t \vec y_{\tau,h}(t_{n,l}^{\text G}), \vec \phi_h \rangle - \langle \vec z_{\tau,h}(t_{n,l}^{\text G}), \vec \phi_h \rangle)\,.
\end{align*}
Thus, we have that 
\begin{equation}
\label{Eq:Prep2_01}
\langle \partial_t \vec y_{\tau,h}(t_{n,l}^{\text G}) - \vec z_{\tau,h}(t_{n,l}^{\text G}), \vec \phi_h \rangle = 0
\end{equation}
for all $\vec \phi_h \in \vec V_h^{r+1}$. Choosing $\vec \phi_h = \partial_t \vec y_{\tau,h}(t_{n,l}^{\text G}) - \vec z_{\tau,h}(t_{n,l}^{\text G})$ in \eqref{Eq:Prep2_01}, proves the assertion \eqref{Eq:Aux02}.
\end{mproof}

\section{Error analysis}
\label{Sec:ErrAna}

Here we derive our error estimate \eqref{Intro:MainEst} for the scheme \eqref{Eq:DPV_0} or \eqref{Eq:DPQ_0}, respectively. In \eqref{Eq:HPS},  let $\vec v:= \partial_t \vec u$. Let $(\vec w_1,\vec w_2)^\top$ be given by Def.~\ref{Def:W}.  We put $\vec U = (\vec u,\vec v)^\top$ and $\vec U_{\tau,h} = (\vec u_{\tau,h},\vec v_{\tau,h})^\top$. We split the error by 
\begin{equation}
\label{Eq:ErrU_00}
\vec U - \vec U_{\tau,h} = \begin{pmatrix}
\vec u - \vec u_{\tau,h} \\ \vec v - \vec v_{\tau,h}  
\end{pmatrix}
= \begin{pmatrix}
\vec u - \vec w_1\\ \vec v - \vec w_2
\end{pmatrix}
+ \begin{pmatrix}
\vec w_1 - \vec u_{\tau,h} \\ \vec w_2- \vec v_{\tau,h}	
\end{pmatrix} =: \begin{pmatrix} \vec \eta_1\\ \vec \eta_2 \end{pmatrix} + \begin{pmatrix} \vec E_{\tau,h}^1\\ \vec E_{\tau,h}^2 \end{pmatrix} = \vec \eta + \vec E_{\tau,h}
\end{equation}
and 
\begin{equation}
\label{Eq:ErrP_00}
p - p_{\tau,h} = p - I_\tau R_h p + I_\tau R_h p - p_{\tau,h} =: \omega +  e_{\tau,h}\,.
\end{equation}
For some quantity $\vec Z = (\vec Z_1,\vec Z_2)^\top\in \vec H^1_0 \times \vec L^2$ we define the norm
\begin{equation}
\label{Eq:DefTN}
||| \vec Z ||| :=  ( \| \nabla \vec Z_1\|^2 + \| \vec Z_2\|^2)^{1/2}
\end{equation}
and the weighted (elastic) energy norm 
\begin{equation}
	\label{Eq:DefTNe}
	||| \vec Z |||_e :=  \Big(\frac{1}{2} \langle \vec C \vec \varepsilon (\vec Z_1), \vec \varepsilon (\vec Z_1) \rangle + \frac{\rho}{2} \langle \vec Z_2 , \vec Z_2 \rangle \Big)^{1/2} \,.
\end{equation}
By Korn's inequality (cf.\ \cite{F47})  along with the positive definiteness of $\vec C$, these norms are equivalent in sense that for $\vec Z = (\vec Z_1,\vec Z_2)^\top\in \vec H^1_0 \times \vec L^2$ there holds that
\begin{equation}
\label{Eq:Nequi}
	c_1 	||| \vec Z ||| \leq ||| \vec Z |||_e \leq c_2  ||| \vec Z ||| 
\end{equation}
with some positive constants $c_1$ and $c_2$. Finally, for some scalar-valued function $q\in L^2$ we define the weighted $L^2$-norm   
\begin{equation}
\label{Eq:DefTNe4s}
||| q |||_e := \Big(\frac{c_0}{2} \langle q,q \rangle \Big)^{1/2}\,.
\end{equation}

We start with providing estimates for the projection errors $\vec \eta$ and $\omega$ of \eqref{Eq:ErrU_00} and \eqref{Eq:ErrP_00}, respectively. For the Lagrange interpolation \eqref{Eq:DefLagIntOp}, $s\in  \{2,\infty\}$ and $m\in \{0,1\}$ we recall that (cf.~\cite{H91})
\begin{equation}
\label{Eq:ErrLI}
\| f - I_\tau f \|_{L^s(I_n;H^m)}   \leq c \tau_n^{k+1} \| \partial_t^{k+1} f\|_{L^s(I_n;H^m)}\,. 	
\end{equation}
For the elliptic projections \eqref{Def:vRh} and \eqref{Def:Rh} onto $V_h^r$ and $\vec V_h^{r+1}$, respectively, we have that (cf., e.g., \cite{BS94})
\begin{subequations}
\label{Eq:EEP}
\begin{alignat}{2}
\label{Eq:Omega}
\| p - R_h p\| + h \| \nabla (p - R_h p) \| & \leq c h^{r+1} \| p \|_{r+1} \,,\\[1ex]
\label{Eq:vmvRhv}
\| \vec v - \vec R_h \vec v \|  + h \| \nabla (\vec v - \vec R_h \vec v) \| & \leq c h^{r+2} \| \vec v \|_{r+2} \,.
\end{alignat}
\end{subequations}

\begin{lem}[Estimates of $\vec \eta$]
\label{Lem:EstEtaOmega}
For $\vec \eta = (\vec u - \vec w_1, \vec v - \vec w_2)^\top$ with $(\vec w_1, \vec w_2)^\top$ of \eqref{Eq:DefW} and $s=2$ or $s=\infty$, there holds that  
\begin{subequations}
\label{Eq:EtauEtav}
\begin{alignat}{2}					
\label{Eq:Etau}
\| \vec u - \vec w_1 \|_{L^s(I_n;\vec L^2)}  & \leq c (\tau_n^{k+1} \mathcal C_{t,s}^{n,1} + h^{r+2} \mathcal C_{\vec x,s}^{n,1}) \,,\\[1ex]
\label{Eq:Etav}
\| \vec v - \vec w_2 \|_{L^s(I_n;\vec L^2)}  & \leq c (\tau_n^{k+1} \mathcal C_{t,s}^{n,2} + h^{r+2}\mathcal C_{\vec x,s}^{n,2}) \,,\\[1ex]
\label{Eq:EtauH1}
\| \vec u - \vec w_1 \|_{L^s(I_n;\vec H^1)}  & \leq c( \tau_n^{k+1} \mathcal C_{t,s}^{n,3} + h^{r+1} \mathcal C_{\vec x,s}^{n,3} )\,,\\[1ex]
\label{Eq:EtavH1}
\| \vec R_h \vec v - \vec w_2 \|_{L^s(I_n;\vec H^1)}  & \leq c (\tau_n^{k+1} \mathcal C_{t,s}^{n,4} + h^{r+1}\mathcal C_{\vec x,s}^{n,4} )\,,
\end{alignat}
\end{subequations}
where the constants in \eqref{Eq:EtauEtav} are given by $\mathcal C_{t,s}^{n,1} :=  \|\partial_t^{k+1} \vec u\|_{L^s(I_n;\vec L^2)}+ \mathcal C_{t,s}^{n,2}$, $\mathcal C_{t,s}^{n,2}  := \|\partial_t^{k+2} \vec u\|_{L^s(I_n;\vec L^2)}$, $\mathcal C_{t,s}^{n,3} :=  \|\partial_t^{k+1} \vec u\|_{L^s(I_n;\vec H^1)}+\tau_n \mathcal C_{t,s}^{n,4}$, $\mathcal C_{t,s}^{n,4} := \|\partial_t^{k+2} \vec u\|_{L^s(I_n;\vec H^1)}$, $\mathcal C_{\vec x,s}^{n,1} := \|\vec u\|_{L^s(I_n;\vec H^{r+2})}+\tau_n\mathcal C_{\vec x,s}^{n,2}$, $\mathcal C_{\vec x,s}^{n,2}  := \|\partial_t \vec u\|_{L^s(I_n;\vec H^{r+2})}+\tau_n\|\partial_t^2 \vec u\|_{L^s(I_n;\vec H^{r+2})}$, $\mathcal C_{\vec x,s}^{n,3} := \|\vec u\|_{L^s(I_n;\vec H^{r+2})} + \tau_n \mathcal C_{\vec x,s}^{n,4}$ and $\mathcal C_{\vec x,s}^{n,4} := \|\partial_t \vec u\|_{L^s(I_n;\vec H^{r+2})}$.
\end{lem}

\begin{mproof}
For scalar-valued functions, estimates \eqref{Eq:Etau} and \eqref{Eq:Etav} are proved in \cite[Lem.~3.3]{KM04} and \eqref{Eq:EtauH1} and \eqref{Eq:EtavH1} in \cite[Appendix]{BKRS20}. The estimates \eqref{Eq:EtauEtav} hold similarly in the vector-valued case of Def.~\ref{Def:W}. 
\end{mproof}

Next, we derive variational equations satisfied by the discretization errors $\vec E_{\tau,h}$ and $e_{\tau,h}$.

\begin{lem}[Variational equations for $ \vec E_{\tau,h}$ and $e_{\tau,h}$]
\label{Lem:EE}
Let 
\begin{equation}
\label{Eq:EE00}
\vec T_I^n : =  I_\tau \int_{t_n-1}^t \partial_t \vec u - I_\tau\partial_t \vec u \ud s\,, \;\;  \vec  T_{II}^n : = \rho \partial_t^2 \vec u - \rho \partial_t \vec w_2\,, \;\;  \vec T_{III}^n :  =  I_\tau \vec u - \vec u\,, 
\; \; \vec T_{IV}^n :  = \rho \vec f - I_\tau (\rho \vec f)\,, \;\; T_{V}^n :  = g - I_\tau g\,,
\end{equation}	
where $I_\tau$ is the Lagrange interpolation operator satisfying \eqref{Eq:DefLagIntOp}. Then, for $n=1,\ldots, N$ the errors $\vec E_{\tau,h}{}_{|I_n}$ and $e_{\tau,h}{}_{|I_n}$ of \eqref{Eq:ErrU_00} and \eqref{Eq:ErrP_00}, respectively, satisfy the equations
\begin{subequations}
	\label{Eq:EE01}
	\begin{alignat}{1}
		\label{Eq:EE02}
		&\int_{I_n}	\llangle \vec D\partial_t \vec E_{\tau,h}, \vec \Phi_{\tau,h} \rrangle + \llangle \vec{\mathcal L_h} \vec E_{\tau,h}, \vec \Phi_{\tau,h} \rrangle \ud t - \alpha \int_{I_n}	 \langle e_{\tau,h}, \nabla \cdot \vec \Phi_{\tau,h}^2 \rangle \ud t =  \int_{I_n}	\langle \vec T_{IV}^n, \Phi_{\tau,h}^2 \rangle \ud t \\[1ex]
		\nonumber
		&\qquad  - \int_{I_n}  \langle \vec A_h \vec T_I^n, \vec \Phi_{\tau,h}^2\rangle \ud t - \int_{I_n}  \langle \vec T_{II}^n, \vec \Phi_{\tau,h}^2\rangle \ud t + \int_{I_n}  \langle \vec A_h \vec T^n_{III}, \vec \Phi_{\tau,h}^2\rangle \ud t  + \alpha \int_{I_n}	\langle \omega, \nabla \cdot \vec \Phi_{\tau,h}^2 \rangle \ud t  \,, \\[2ex]
		\label{Eq:EE03}
		& \int_{I_n} c_0 \langle \partial_t e_{\tau,h}, \psi_{\tau,h}\rangle+ \langle B_h e_{\tau,h}, \psi_{\tau,h}\rangle \ud t  +  \alpha \int_{I_n} \langle \nabla \cdot \partial_t \vec E_{\tau,h}^1, \psi_{\tau,h}\rangle \ud t \\[1ex]
		\nonumber
		& \quad = \int_{I_n} \langle T_{V}^n, \psi_{\tau,h}\rangle - \int_{I_n} c_0 \langle \partial_t \omega, \psi_{\tau,h}\rangle \ud t  - \alpha \int_{I_n}\langle \nabla \cdot \partial_t \vec \eta_1, \psi_{\tau,h}\rangle \ud t - \int_{I_n} \langle \vec K \nabla (p-I_\tau p) , \nabla \psi_{\tau,h}\rangle \ud t
	\end{alignat}
\end{subequations}	 
for all $ \vec \Phi_{\tau,h}\in (\mathbb P_{k-1}(I_n;V_{h}^{r+1}))^d\times (\mathbb P_{k-1}(I_n;V_{h}^{r+1}))^d$ and ${\psi_{\tau,h} \in \mathbb P_{k-1}(I_n;V_{h}^{r})}$. 
\end{lem}

\begin{mproof}
Let $\vec v= \partial_t \vec u$. Rewriting \eqref{Eq:HPS_1} as a first-order in time system, substracting \eqref{Eq:DPV_0} from the weak form of the resulting first-order in time, continuous system and using the splitting \eqref{Eq:ErrU_00} and \eqref{Eq:ErrP_00} of the errors we get that 
\begin{subequations}
	\label{Eq:EE10}
	\begin{alignat}{1}		
		\label{Eq:EE11}
		&\int_{I_n}	\llangle \vec D \partial_t \vec E_{\tau,h}, \vec \Phi_{\tau,h} \rrangle + \llangle \vec{\mathcal L_h} \vec E_{\tau,h}, \vec \Phi_{\tau,h} \rrangle - \alpha \langle e_{\tau,h}, \nabla \cdot \vec \Phi_{\tau,h}^2 \rangle \ud t\\
		\nonumber
		&=  \int_{I_n}	\llangle \vec F - I_\tau \vec F, \Phi_{\tau,h} \rangle \ud t - \int_{I_n} \llangle \vec D \partial_t \vec \eta, \vec \Phi_{\tau,h} \rrangle + \llangle \vec{\mathcal L_h} \vec \eta, \vec \Phi_{\tau,h} \rrangle - \alpha \langle \omega, \nabla \cdot \vec \Phi_{\tau,h}^2 \rangle \ud t\,, \\[1ex]
		\label{Eq:EE12}
		& \int_{I_n} c_0 \langle \partial_t e_{\tau,h}, \psi_{\tau,h}\rangle + \alpha \langle \nabla \cdot \partial_t \vec E_{\tau,h}^1, \psi_{\tau,h}\rangle + \langle B_h e_{\tau,h}, \psi_{\tau,h}\rangle \ud t \\
		\nonumber
		& = \int_{I_n} \langle g - I_\tau g, \psi_{\tau,h}\rangle - \int_{I_n} c_0 \langle \partial_t \omega, \psi_{\tau,h}\rangle + \alpha \langle \nabla \cdot \partial_t \vec \eta_1, \psi_{\tau,h}\rangle + \langle B_h \omega , \psi_{\tau,h}\rangle \ud t
	\end{alignat}
\end{subequations}	 
for all $ \vec \Phi_{\tau,h}\in (\mathbb P_{k-1}(I_n;V_{h}^{r+1}))^d\times (\mathbb P_{k-1}(I_n;V_{h}^{r+1}))^d$ and ${\psi_{\tau,h} \in \mathbb P_{k-1}(I_n;V_{h}^{r})}$. 

Next, we rewrite some of the terms in  \eqref{Eq:EE10}. Firstly, from \eqref{Eq:ErrU_00} along with \eqref{Eq:DefLh} we find that 
\begin{equation}
\begin{aligned}
\label{Eq:EE13}
\int_{I_n} \llangle \vec D \partial_t \vec \eta, \vec \Phi_{\tau,h} \rrangle + \llangle \vec{\mathcal L_h} \vec \eta, \vec \Phi_{\tau,h} \rrangle  \ud t &  = \int_{I_n} \langle \partial_t \vec u - \partial_t \vec w_1 - \vec v + \vec w_2, \vec \Phi_{\tau,h}^1 \rangle \ud t\\[1ex]
&\quad  + \int_{I_n}  \langle \rho \partial_t \vec v - \rho \partial_t \vec w_2  + \vec A_h (\vec u - \vec w_1), \vec \Phi_{\tau,h}^2\rangle \ud t\,.
\end{aligned}
\end{equation}
Recalling that $\vec v = \partial_t \vec u$ and Lem.~\ref{Lem:Relwj}, we get for the first term on the right-hand side of \eqref{Eq:EE13} that 
\begin{equation}
\label{Eq:EE14}
\int_{I_n} \langle \partial_t \vec u - \partial_t \vec w_1 - \vec v + \vec w_2, \vec \Phi_{\tau,h}^1 \rangle \ud t = 0 
\end{equation}
for all $ \vec \Phi_{\tau,h}^1\in (\mathbb P_{k-1}(I_n;V_{h}^{r+1}))^d$. Let now 
\begin{equation}
\label{Eq:EE15}
\vec z(t) := \int_{t_{n-1}}^t \vec w_2(s)\ud s + \vec R_h \vec u(t_{n-1})\,.
\end{equation}
Then, by definition we have that 
\begin{equation}
\label{Eq:EE16}
\vec w_1{}_{|I_n} = I_\tau \vec z\,. 
\end{equation}
For the last term on the right-hand side of \eqref{Eq:EE13} we get by \eqref{Eq:EE15}, \eqref{Eq:EE16} and \eqref{Eq:DefW} along with \eqref{Def:vRh} and \eqref{Eq:DefAh} that  
\begin{equation}
\begin{aligned}
\label{Eq:EE17}
& \int_{I_n}  \langle \vec A_h \vec w_1, \vec \Phi_{\tau,h}^2\rangle \ud t\\[1ex] 
& = \frac{\tau_n}{2}\sum_{\mu=0}^k \hat \omega_\mu  \langle \vec A_h \vec z(t_{n,\mu}), \vec \Phi_{\tau,h}^2(t_{n,\mu})\rangle 
 = \frac{\tau_n}{2}\sum_{\mu=0}^k \hat \omega_\mu  \Big\langle \vec A_h \Big(\int_{t_{n-1}}^{t_{n,\mu}} I_\tau \vec R_h \partial_t \vec u \ud s +  \vec u(t_{n-1})\Big), \vec \Phi_{\tau,h}^2(t_{n,\mu})\Big\rangle\\[1ex]
& = - \int_{I_n}  \Big\langle \vec A_h \vec T^n_I, \vec \Phi_{\tau,h}^2\Big\rangle \ud t + \int_{I_n}  \langle \vec A_h \vec u, \vec \Phi_{\tau,h}^2\rangle \ud t + \int_{I_n}  \langle \vec A_h \vec T^n_{III}, \vec \Phi_{\tau,h}^2\rangle \ud t
\end{aligned}
\end{equation}
with $\vec T^n_I$ and $\vec T^n_{III}$ being defined in \eqref{Eq:EE00}. Combining now \eqref{Eq:EE13} with \eqref{Eq:EE14} and \eqref{Eq:EE17} yields that
\begin{equation}
\label{Eq:EE18}
\begin{aligned}
\int_{I_n} \llangle \vec D \partial_t \vec \eta, \vec \Phi_{\tau,h} \rrangle + \llangle \vec{\mathcal L_h} \vec \eta, \vec \Phi_{\tau,h} \rrangle  \ud t & = \int_{I_n}  \langle \vec A_h \vec T_I^n, \vec \Phi_{\tau,h}^2\rangle \ud t + \int_{I_n}  \langle \vec T_{II}^n, \vec \Phi_{\tau,h}^2\rangle \ud t - \int_{I_n}  \langle \vec A_h \vec T^n_{III}, \vec \Phi_{\tau,h}^2\rangle \ud t 
\end{aligned}
\end{equation}
with $\vec T^n_{II}$ being defined in \eqref{Eq:EE00}. Together, \eqref{Eq:EE11} and \eqref{Eq:EE18} prove the assertion \eqref{Eq:EE02}. 

For the last of the terms on the right-hand side of \eqref{Eq:EE03} it holds by \eqref{Def:Rh} that 
\begin{equation}
	\label{Eq:EE20}
	\begin{aligned}
	\int_{I_n} \langle B_h \omega , \psi_{\tau,h}\rangle \ud t & = \int_{I_n} \langle \vec K \nabla (p-I_\tau R_h p) , \nabla \psi_{\tau,h}\rangle \ud t
	=  \int_{I_n} \langle \vec K \nabla (p-I_\tau p) , \nabla \psi_{\tau,h}\rangle \ud t\\[1ex]& \quad  +  \int_{I_n} \langle \vec K \nabla (I_\tau p- R_h I_\tau p) , \nabla \psi_{\tau,h}\rangle \ud t = \int_{I_n} \langle \vec K \nabla (p-I_\tau p) , \nabla \psi_{\tau,h}\rangle \ud t 
	\end{aligned}
\end{equation}
for all $\psi_{\tau,h}\in \mathbb P_{k-1}(I_n;V_h^d)$. Together, \eqref{Eq:EE12} and \eqref{Eq:EE20} prove the assertion \eqref{Eq:EE03}.
\end{mproof}

The following lemma provides estimates for the terms $\vec T_I$, $\vec T_{II}^n$, $\vec T_{III}^n$ of  \eqref{Eq:EE00} and $\partial_t \vec \eta_1$ of \eqref{Eq:EE01}.

\begin{lem}[Estimation of $\vec T_I$, $\vec T_{II}^n$, $\vec T_{III}^n$ and $\partial_t \vec \eta_1$]
\label{Lem:EstTni}
For $\vec T_I^n$, $\vec T_{II}^n$, $\vec T_{III}^n$ and $\partial_t \vec \eta_1$ there holds that 
\begin{subequations}
\label{Eq:EstTni0}
\begin{alignat}{3}
\label{Eq:EstTni1}
\|\vec A \vec T_{I}^n\|_{L^2(I_n;\vec L^2)} &  \leq c \tau_n^{k+1} \| \vec A \partial_t^{k+1} \vec u \|_{L^2(I_n;\vec L^2)}\,,
\\[1ex]
\label{Eq:EstTni2}
\left| \int_{I_n} \langle \vec T_{II}^n, \vec \phi_{\tau,h}\rangle \ud t \right|  & \leq c \big(\tau_n^{k+1}\|\partial_t^{k+3} \vec u\|_{L^2(I_n;\vec L^2)}  + h^{r+2}\|\partial_t^2 \vec u\|_{L^2(I_n;\vec H^{r+2})} \big)\|\vec \phi_{\tau,h}\|_{L^2(I_n;\vec L^2)}\,,\\[1ex]
\label{Eq:EstTni3}
\|\vec A_h \vec T_{III}^n\|_{L^2(I_n;\vec L^2)} & \leq c  \tau_n^{k+1}\| \vec A \partial_t^{k+1} \vec u \|_{L^2(I_n;\vec L^2)}\,,\\[1ex]
\label{Eq:Estdiveta}
 \left| \int_{I_n} \langle \nabla \cdot \partial_t \vec \eta_1, \psi_{\tau,h}\rangle \ud t\right| & \leq c \Big(\tau_n^{k+1}\|\partial_t^{k+2} \vec u\|_{L^2(I_n;\vec H^1)} + h^{r+1}  \|\partial_t \vec u\|_{L^2(I_n;\vec H^{r+2})}\Big)\| \psi_{\tau,h}\|_{L^2(I_n;L^2)}\,,\\[1ex]
\label{Eq:Estdtomega}
\left| \int_{I_n} \langle  \partial_t \omega, \psi_{\tau,h}\rangle \ud t\right| & \leq c \Big(\tau_n^{k+1}\|\partial_t^{k+2} p \|_{L^2(I_n;L^2)} + h^{r+1}  \|\partial_t p \|_{L^2(I_n;H^{r+1})}\Big)\| \psi_{\tau,h}\|_{L^2(I_n;L^2)}
\end{alignat}
for $\vec \phi_{\tau,h} \in (\mathbb P_{k-1}(I_n;V^{r+1}_h))^d$ in \eqref{Eq:EstTni2} and $\psi_{\tau,h} \in \mathbb P_{k-1}(I_n;V^{r}_h)$ in \eqref{Eq:Estdiveta}. 
\end{subequations}
\end{lem}

\begin{mproof}
The inequalities \eqref{Eq:EstTni1} to \eqref{Eq:EstTni3} can be proved along the lines of \cite[Lem.~3.3, Eqs.\ (3.12) to (3.14)]{KM04} that are shown for scalar-valued functions. It remains to prove \eqref{Eq:Estdiveta} for $\vec \eta_1 = \vec u -\vec w_1$ and \eqref{Eq:Estdtomega} for $\omega = p - I_\tau R_hp$. From the first of the definitions in \eqref{Eq:DefW} it follows that 
\begin{equation}
\label{Eq:Estdiveta01}
\vec \eta_1 = \vec u - \vec w_1 = \vec u - I_\tau \vec u + I_\tau \vec u - I_\tau (\vec R_h \vec u) - I_\tau 	\int_{t_{n-1}}^t (\vec w_2 - \partial_t \vec R_h \vec u) \ud s\,.
\end{equation}
By \eqref{Eq:Estdiveta01} we then get that 
\begin{equation}
\label{Eq:Estdiveta02}
\begin{aligned}
& \int_{I_n} \langle \nabla \cdot \partial_t \vec \eta_1, \psi_{\tau,h}\rangle \ud t =  \int_{I_n} \langle \nabla \cdot \partial_t (\vec u - I_\tau \vec u ), \psi_{\tau,h}\rangle \ud t+ \int_{I_n} \langle \nabla \cdot \partial_t I_\tau (\vec u  - \vec R_h \vec u ), \psi_{\tau,h}\rangle \ud t\\[1ex]
& \qquad  + \int_{I_n} \Big\langle \nabla \cdot \partial_t I_\tau 	\int_{t_{n-1}}^t (\vec w_2 - \partial_t \vec R_h \vec u) \ud s, \psi_{\tau,h}\Big\rangle \ud t   =: \Gamma_1 + \Gamma_2 + \Gamma_3\,.
 \end{aligned}
\end{equation}

We start with estimating $\Gamma_1$. Firstly, let $k\geq 2$. Using integration by parts in time and recalling that the endpoints of $I_n$ are included in the set of Gauss--Lobatto quadrature points  of $I_n$, we get that 
\begin{equation*}
\Gamma_1 = \int_{I_n} \langle \nabla \cdot \partial_t (\vec u - I_\tau \vec u ), \psi_{\tau,h}\rangle \ud t = - \int_{I_n} \langle \nabla \cdot (\vec u - I_\tau \vec u), \partial_t \psi_{\tau,h}\rangle \ud t\,.
\end{equation*}
Let now $I_\tau^{k+1}$ denote the Lagrange interpolation operator at the $k+2$ points of $\overline I_n = [t_{n-1},t_n]$ consisting of the $k+1$ Gauss--Lobatto quadrature nodes $t_{n,\mu}$, for $\mu=0,\ldots,k$, and a further node in $(t_{n-1},t_n)$ that is distinct from the previous ones. Then, $(I_\tau^{k+1}\vec u)\partial_t \psi_{\tau,h}$ is a polynomial of degree $2k-1$ in $t$, such that 
\begin{equation*}
\int_{I_n} \langle \nabla \cdot (\vec u - I_\tau \vec u), \partial_t \psi_{\tau,h}\rangle \ud t = \int_{I_n} \langle \nabla \cdot (\vec u - I_\tau^{k+1} \vec u), \partial_t \psi_{\tau,h}\rangle \ud t\,.
\end{equation*}
Using integration by parts, the stability of the operator $I_\tau^{k+1}$ in the norm of $L^2(I_n;H^1)$, we have that
\begin{equation}
\label{Eq:Estdiveta03}
\begin{aligned}
|\Gamma_1| & \leq \left|\int_{I_n} \langle \nabla \cdot \partial_t (\vec u - I_\tau^{k+1} \vec u), \psi_{\tau,h}\rangle \ud t\right|\\[1ex] 
& \leq \| \partial_t (\vec u - I_\tau^{k+1} \vec u) \|_{L^2(I_n;\vec H^1)} \|\psi_{\tau,h}\|_{L^2(I_n;L^2)}\\[1ex] 
&  \leq c \tau_n^{k+1}\| \partial_t^{k+2} \vec u\|_{L^2(I_n;\vec H^1)} \|\psi_{\tau,h}\|_{L^2(I_n;L^2)}\,.
\end{aligned}
\end{equation}
For $k=1$, we have that $\partial_t I_\tau \vec u, \psi_{\tau,h}\in \mathbb P_0(I_n;V_h^r)$ with $\partial_t I_\tau \vec u = (\vec u(t_n)-\vec u(t_{n-1}))/\tau_n$. It follows that 
\begin{equation}
	\label{Eq:Estdiveta03k1}
	\begin{aligned}
\Gamma_1 = \Big\langle \nabla \cdot \int_{I_n}  (\partial_t \vec u - \partial_t I_\tau\vec u)\ud t , \psi_{\tau,h}\Big\rangle
=  \langle \nabla \cdot (\vec u(t_n)-\vec u(t_{n-1})-(\vec u(t_n)-\vec u(t_{n-1}))) , \psi_{\tau,h}\Big\rangle = 0\,.
\end{aligned}
\end{equation}

Next, we estimate $\Gamma_2$. For this we introduce the abbreviation $\vec \xi := \vec u -\vec R_h \vec u$. The Lagrange interpolant $I_\tau$ satisfies the stability results (cf.~\cite[Eqs.~(3.15) and (3.16)]{KM04})
\begin{subequations}
\label{Eq:Estdiveta035}
\begin{alignat}{2}
\label{Eq:Estdiveta0351}
\|  I_\tau w \|_{L^2(I_n;L^2)} & \leq c \| w \|_{L^2(I_n;L^2)}  + c \tau_n \|  \partial_t w \|_{L^2(I_n;L^2)} \,, \\[1ex]
\label{Eq:Estdiveta0352}
 \Big\|  \int_{t_{n-1}}^t w \ud s \Big \|_{L^2(I_n;L^2)}  & \leq c \tau_n  \|  w \|_{L^2(I_n;L^2)} \,.
\end{alignat}
\end{subequations}
By the $H^1$--$L^2$ inverse inequality $\| w'\|_{L^2(I_n;\R)}\leq c \tau_n^{-1}\| w\|_{L^2(I_n;\R)}$, the stability results \eqref{Eq:Estdiveta035}, the error estimate \eqref{Eq:vmvRhv} and viewing $\vec \xi(t_{n-1}^+)$ as a function constant in time we find that 
\begin{equation}
\label{Eq:Estdiveta04}
\begin{aligned}
|\Gamma_2 | & = \Big| \int_{I_n} \langle \nabla \cdot \partial_t I_\tau \vec \xi, \psi_{\tau,h}\rangle \ud t\Big| = \Big| \int_{I_n} \langle \nabla \cdot \partial_t I_\tau (\vec \xi- \vec \xi(t_{n-1}^+)), \psi_{\tau,h}\rangle \ud t\Big|\\[1ex]
&  = \Big| \int_{I_n} \langle \nabla \cdot \partial_t I_\tau \int_{t_{n-1}}^t \partial_t \vec \xi\ud s, \psi_{\tau,h}\rangle \ud t\Big|\\[1ex] 
& \leq c \tau_n^{-1} \Big\| I_\tau \int_{t_{n-1}}^{t} \nabla \cdot \partial_t \vec \xi \ud s
 \Big\|_{L^2(I_n;L^2)} \|\psi_{\tau,h}\|_{L^2(I_n;L^2)}\\[1ex]
 &
\leq c h^{r+1}\| \partial_t \vec u \|_{L^2(I_n;\vec H^{r+2})}\| \psi_{\tau,h}\|_{L^2(I_n;L^2)}\,.
\end{aligned}
\end{equation}
Finally, we estimate $\Gamma_3$. By the arguments of \eqref{Eq:Estdiveta04} it follows for $\Gamma_3$ that 
\begin{equation*}
\begin{aligned}
|\Gamma_3 | & = \Big| \int_{I_n} \Big\langle \nabla \cdot \partial_t I_\tau \int_{t_{n-1}}^t (\vec w_2 - \partial_t \vec R_h \vec u) \ud s, \psi_{\tau,h}\Big\rangle \ud t\Big| \leq c \|\vec w_2  - \vec R_h (\partial_t \vec u)  \|_{L^2(I_n;\vec H^1)}\|\psi_{\tau,h}\|_{L^2(I_n;L^2)}\,.
\end{aligned}
\end{equation*}
Employing \eqref{Eq:EtavH1} with $\vec v =\partial_t \vec u$, we obtain that
\begin{equation}
\label{Eq:Estdiveta05}
\begin{aligned}
|\Gamma_3 | & \leq c \Big(\tau_n^{k+1}\|\partial_t^{k+2} \vec u \|_{L^2(I_n;\vec H^1)} + h^{r+1}\| \partial_t \vec u \|_{L^2(I_n;\vec H^{r+2})}\Big) \|\psi_{\tau,h}\|_{L^2(I_n;L^2)}\,.
\end{aligned}
\end{equation} 
Now, combining \eqref{Eq:Estdiveta02} with \eqref{Eq:Estdiveta03}, \eqref{Eq:Estdiveta03k1}, \eqref{Eq:Estdiveta04} and \eqref{Eq:Estdiveta05} proves the assertion \eqref{Eq:Estdiveta}.  Estimate \eqref{Eq:Estdtomega} can be shown similarly to \eqref{Eq:Estdiveta} along the lines of \eqref{Eq:Estdiveta01} to \eqref{Eq:Estdiveta04}.
\end{mproof}

Next, we prove a stability estimate for the error $||| \vec E_{\tau,h} (t_n)|||_e^2 + |||e_{\tau,h}(t_n)|||_e^2 $.

\begin{lem}[Stability estimate]
\label{Lem:ES}
Let $n=1,\ldots, N$ and  
\begin{equation}
	\label{Eq:ES00}
\delta_n := \alpha \langle \omega(t_n), \nabla \cdot \vec E_{\tau,h}^1(t_n) \rangle \quad \text{and} \quad \delta_{n-1}^+:= \alpha \langle \omega(t_{n-1}^+), \nabla \cdot \vec E_{\tau,h}^1(t_{n-1}^+) \rangle\,,
\end{equation}
where the errors $\vec E_{\tau,h}$, $e_{\tau,h}$ and $\omega$ are defined in \eqref{Eq:ErrU_00} and \eqref{Eq:ErrP_00}, respectively. Then, there holds that 
\begin{equation}
\label{Eq:ES01}
\begin{aligned}
||| \vec E_{\tau,h} (t_n)|||_e^2 + ||| e_{\tau,h}(t_n)|||_e^2   & \leq ||| \vec E_{\tau,h} (t_{n-1}^+) |||_e^2 + |||e_{\tau,h}(t_{n-1}^+)|||_e^2 +\delta_n - \delta_{n-1}^+ + c ||| \vec E_{\tau,h}|||^2_{L^2(I_n;\vec L^2)} + c \|e_{\tau,h}\|^2_{L^2(I_n;L^2)} \\[1ex]
	& \qquad + c \tau_n^{2(k+1)} (\mathcal E_t^{n,1})^2 + c h^{2(r+1)} (\mathcal E_{\vec x}^{n,1})^2 + c h^{2(r+2)} (\mathcal E_{\vec x}^{n,2})^2 	
\end{aligned}
\end{equation}
with $\mathcal E_t^{n,1}:=\mathcal E_{\vec u}^{I,n}+\mathcal E_{\vec u,t}^{II,n}+\mathcal E_{\vec u}^{III,n}+\mathcal E_{\vec f}^n+\mathcal E_{g}^n + \mathcal E_{\omega,t}^{n} + \mathcal E_{\vec \eta,t}^{n}+\mathcal E_{p,t}^n$, $\mathcal E_{\vec x}^{n,1} := \mathcal E_{\omega,\vec x}^n + \mathcal E_{\vec \eta,\vec x}^{n}$, $\mathcal E_{\vec x}^{n,2} := \mathcal E_{\vec u,\vec x}^{II,n}$, where 
\begin{equation*}
\label{Eq:DefMcE}
\begin{aligned}
\mathcal E_{\vec u}^{I,n} & :=\| \partial_t^{k+1} \vec u\|_{L^2(I_n;\vec H^2)}\,, & \mathcal E_{\vec u,t}^{II,n} & :=\| \partial_t^{k+3} \vec u\|_{L^2(I_n;\vec L^2)}\,, & \mathcal E_{\vec u,\vec x}^{II,n} & :=\| \partial_t^{2} \vec u\|_{L^2(I_n;\vec H^{r+2})}\,, & 
\mathcal E_{\vec u}^{III,n} & :=\| \partial_t^{k+1} \vec u\|_{L^2(I_n;\vec H^2)}\,, \\[1ex]
 \mathcal E_{\vec f}^n & :=\| \partial_t^{k+1} \vec f\|_{L^2(I_n;\vec L^2)}\,, & \mathcal E_{g}^n & : =\| \partial_t^{k+1} g \|_{L^2(I_n;L^2)}\,,  & \mathcal E_{\omega,t}^n & :=\| \partial_t^{k+2} p \|_{L^2(I_n;L^2)}\,, & \mathcal E_{\omega,\vec x}^n & :=\| \partial_t p \|_{L^2(I_n;H^{r+1})}\,, \\[1ex]
\mathcal E_{p,t}^n & := \| \partial_t^{k+1}  p\|_{L^2(I_n;H^2)}\,, & \mathcal E_{\vec \eta,t}^{n} & := \| \partial_t^{k+2} \vec u\|_{L^2(I_n;\vec H^1)}\,, & 
\mathcal E_{\vec \eta,\vec x}^{n} & := \| \partial_t \vec u\|_{L^2(I_n;\vec H^{r+2})} \,.
\end{aligned}
\end{equation*}
\end{lem}

\begin{mproof}
In \eqref{Eq:EE01}, we choose the test functions 
\begin{equation}
	\label{Eq:ES10}
	\vec \Phi_{\tau,h} = \begin{pmatrix} \Pi^{k-1}_\tau & 0 \\ 0 & \Pi^{k-1}_\tau\end{pmatrix}\begin{pmatrix} \vec A_h \vec E_{\tau,h}^1\\ \vec E_{\tau,h}^2 \end{pmatrix} \quad \text{and} \quad \psi_{\tau,h} = \Pi^{k-1}_\tau e_{\tau,h}\,.
\end{equation}
Firstly, we address some of the terms in \eqref{Eq:EE02} for the test function $\vec \Phi_{\tau,h}$ of \eqref{Eq:ES10}. By the exactness of the Gauss quadrature formula \eqref{Eq:GF} for all polynomials in $\mathbb P_{2k-1}(I_n;\R)$ and Lem.~\ref{Lem:PropGF} we deduce that 
\begin{equation}
\label{Eq:ES11}
\begin{aligned}
    & \int_{I_n}	\left\llangle \begin{pmatrix} \vec I_d & \vec 0\\ \vec 0 & \rho\vec I_d \end{pmatrix}\begin{pmatrix} \partial_t \vec E_{\tau,h}^1\\ \partial_t \vec E_{\tau,h}^2\end{pmatrix},  \begin{pmatrix}\Pi^{k-1}_\tau \vec A_h \vec E_{\tau,h}^1\\  \Pi^{k-1}_\tau\vec E_{\tau,h}^2 \end{pmatrix}\right\rrangle \ud t  =  \frac{\tau_n}{2}\sum_{\mu = 1}^k \hat \omega_\mu^{\text G} \left\llangle \begin{pmatrix} \partial_t \vec E_{\tau,h}^1(t_{n,\mu}^{\text G})\\ \partial_t \vec E_{\tau,h}^2(t_{n,\mu}^{\text G})\end{pmatrix},  \begin{pmatrix} \vec A_h \vec E_{\tau,h}^1(t_{n,\mu}^{\text G})\\  \vec E_{\tau,h}^2(t_{n,\mu}^{\text G}) \end{pmatrix}\right\rrangle\\[1ex]
& = \int_{I_n} \underbrace{\langle \partial_t \vec E_{\tau,h}^1,  \vec A_h \vec E_{\tau,h}^1\rangle}_{= \frac{1}{2} \frac{d}{dt} \langle \vec A_h \vec E_{\tau,h}^1,  \vec E_{\tau,h}^1\rangle } + \rho \underbrace{\langle \partial_t \vec E_{\tau,h}^2, \vec E_{\tau,h}^2\rangle}_{= \frac{1}{2} \frac{d}{dt} \langle \vec E_{\tau,h}^2, \vec E_{\tau,h}^2\rangle}\ud t = ||| \vec E_{\tau,h} (t_n)|||_e^2 - ||| \vec E_{\tau,h} (t_{n-1}^+)|||_e^2 \,.
\end{aligned}
\end{equation}
Further, by \eqref{Eq:DefLh}, the exactness of the Gauss quadrature formula \eqref{Eq:GF} for all polynomials in $\mathbb P_{2k-1}(I_n;\R)$, Lem.~\ref{Lem:PropGF} and the symmetry of $\vec A_h$ we have that 
\begin{equation}
\label{Eq:ES12}
\begin{aligned}
& \int_{I_n}\left\llangle \underbrace{\begin{pmatrix} - \vec E_{\tau,h}^2\\ \vec A_h \vec E_{\tau,h}^1\end{pmatrix}}_{=\vec{\mathcal L_h}\vec E_{\tau,h}},  \begin{pmatrix}\Pi^{k-1}_\tau \vec A_h \vec E_{\tau,h}^1\\  \Pi^{k-1}_\tau\vec E_{\tau,h}^2 \end{pmatrix}\right\rrangle \ud t =  \frac{\tau_n}{2}\sum_{\mu = 1}^k \hat \omega_\mu^{\text G} \left\llangle \begin{pmatrix} - \vec E_{\tau,h}^2(t_{n,\mu}^{\text G})\\ \vec A_h \vec E_{\tau,h}^1(t_{n,\mu}^{\text G})\end{pmatrix},  \begin{pmatrix} \vec A_h \vec E_{\tau,h}^1(t_{n,\mu}^{\text G})\\  \vec E_{\tau,h}^2(t_{n,\mu}^{\text G}) \end{pmatrix}\right\rrangle = 0\,.
\end{aligned}
\end{equation}

Next, we recall the definition of the error $(\vec E^1_{\tau,h},\vec E^2_{\tau,h})$ in \eqref{Eq:ErrU_00}. The pair $(\vec w_1,\vec w_2)$ satisfies \eqref{Eq:Relwj} and $(\vec u_{\tau,h},\vec v_{\tau,h})$ fulfills the first of the identities in \eqref{Eq:DPV_1} or \eqref{Eq:DPQ_1}, respectively. Therefore, Lem.~\ref{Lem:upev} can be applied to $(\vec w_1,\vec w_2)$ and $(\vec u_{\tau,h},\vec v_{\tau,h})$ and the conclusion \eqref{Eq:Aux02} holds for both tuples of functions. This implies that 
\begin{equation} 
\label{Eq:ES13}
\vec E_{\tau,h}^2(t_{n,\mu}^{\text G}) = \vec w_2 (t_{n,\mu}^{\text G}) - \vec v_{\tau,h}(t_{n,\mu}^{\text G}) = \partial_t \vec w_1 (t_{n,\mu}^{\text G}) - \partial_t \vec u_{\tau,h}(t_{n,\mu}^{\text G}) = \partial_t \vec E_{\tau,h}^1(t_{n,\mu}^{\text G}) 
\end{equation}
for $\mu=1,\ldots,k$. Using this along with \eqref{Eq:PropGF_0}, it follows that 
\begin{equation}
\label{Eq:ES14}
\begin{aligned}
& \int_{I_n}  \underbrace{\langle e_{\tau,h}, \nabla \cdot \Pi^{k-1}_\tau \vec E_{\tau,h}^2 \rangle}_{ \in \mathbb P_{2k-1}(I_n;\R)} \ud t 
= \frac{\tau_n}{2}\sum_{\mu = 1}^k \hat \omega_\mu^{\text G} \langle e_{\tau,h}(t_{n,\mu}^{\text G}) ,\nabla \cdot \vec E_{\tau,h}^2(t_{n,\mu}^{\text G})  \rangle  \\[1ex]
& \qquad = \frac{\tau_n}{2}\sum_{\mu = 1}^k \hat \omega_\mu^{\text G} \langle e_{\tau,h}(t_{n,\mu}^{\text G}) ,\nabla \cdot \partial_t \vec E_{\tau,h}^1(t_{n,\mu}^{\text G})  \rangle = \int_{I_n} \underbrace{\langle e_{\tau,h}, \nabla \cdot \partial_t \vec E_{\tau,h}^1 \ud t\rangle}_{ \in \mathbb P_{2k-1}(I_n;\R)}  \ud t \,.
\end{aligned}
\end{equation}
By the same arguments and using that $\partial_t \vec E_{\tau,h}^1\in \mathbb (P_{k-1}(I_n;V_h^{r+1}))^d$, we have that 
\begin{equation}
\label{Eq:ES15}
\begin{aligned}
& \int_{I_n} \langle\omega, \nabla \cdot \Pi^{k-1}_\tau \vec E_{\tau,h}^2 \rangle\ud t =  \int_{I_n} \langle \Pi_\tau^{k-1}\omega, \nabla \cdot \Pi^{k-1}_\tau \vec E_{\tau,h}^2 \rangle\ud t  = \frac{\tau_n}{2}\sum_{\mu = 1}^k \hat \omega_\mu^{\text G} \langle  \Pi^{k-1}_\tau \omega(t_{n,\mu}^{\text G}) ,\nabla \cdot \vec E_{\tau,h}^2(t_{n,\mu}^{\text G})  \rangle\\[1ex]  
& = \frac{\tau_n}{2}\sum_{\mu = 1}^k \hat \omega_\mu^{\text G} \langle \Pi^{k-1}_\tau  \omega (t_{n,\mu}^{\text G}) ,\nabla \cdot \partial_t \vec E_{\tau,h}^1(t_{n,\mu}^{\text G})  \rangle = \int_{I_n} \langle \Pi_\tau^{k-1}\omega, \nabla \cdot \partial_t \vec E_{\tau,h}^1 \rangle \ud t  = \int_{I_n} \langle \omega, \nabla \cdot \partial_t \vec E_{\tau,h}^1 \rangle \ud t \,.
\end{aligned}
\end{equation}
Applying integration by parts (for the time variable) to the last term in \eqref{Eq:ES15}, we get that 
\begin{equation}
\label{Eq:ES16}
\begin{aligned}
\int_{I_n}\langle\omega, \nabla \cdot \Pi^{k-1}_\tau \vec E_{\tau,h}^2 \rangle \ud t &  = - \int_{I_n}\langle\partial_t \omega, \nabla \cdot \vec E_{\tau,h}^1 \rangle \ud t + \langle \omega(t_n), \nabla \cdot \vec E_{\tau,h}^1(t_n) \rangle - \langle \omega(t_{n-1}^+), \nabla \cdot \vec E_{\tau,h}^1(t_{n-1}^+) \rangle\,. 
\end{aligned}
\end{equation}

Secondly, we address  some of the terms in \eqref{Eq:EE03} for $\psi_{\tau,h}$  being given by \eqref{Eq:ES10}. Similarly to \eqref{Eq:ES11}, we get that 
\begin{equation}
\label{Eq:ES17}
\begin{aligned}
c_0 \int_{I_n} \langle \partial_t e_{\tau,h}, \Pi_\tau^{k-1} e_{\tau,h} \rangle \ud t & = c_0 \frac{\tau_n}{2} \sum_{\mu=1}^k \langle \partial_t e_{\tau,h}(t_{n,\mu}^{\text G}),  e_{\tau,h}(t_{n,\mu}^{\text G}) \rangle \\[1ex] & = c_0 \int_{I_n} \langle \partial_t e_{\tau,h}, e_{\tau,h} \rangle \ud t \\[1ex] & =  |||e_{\tau,h}(t_n)|||_e^2 - ||| e_{\tau,h}(t_{n-1}^+)|||_e^2\,. 
\end{aligned}
\end{equation}
Further, it holds that 
\begin{equation}
\label{Eq:ES18}
\begin{aligned}
\int_{I_n} \langle B_h e_{\tau,h}, \Pi_\tau^{k-1} e_{\tau,h} \rangle \ud t & = \int_{I_n} \langle B_h \Pi_\tau^{k-1}  e_{\tau,h}, \Pi_\tau^{k-1} e_{\tau,h} \rangle \ud t \,.
\end{aligned}
\end{equation}

Now, adding the equations \eqref{Eq:EE02} and \eqref{Eq:EE03} for the test functions \eqref{Eq:ES10} and using \eqref{Eq:ES11} to \eqref{Eq:ES18} we obtain that 
\begin{equation}
\label{Eq:ES20}
\begin{aligned}
& ||| \vec E_{\tau,h} (t_n)|||_e^2 + ||| e_{\tau,h}(t_n)|||_e^2 + \int_{I_n} \langle B_h \Pi_\tau^{k-1}  e_{\tau,h}, \Pi_\tau^{k-1} e_{\tau,h} \rangle \ud t  =  ||| \vec E_{\tau,h}^1 (t_{n-1}^+)|||_e^2 + ||| e_{\tau,h}(t_{n-1}^+)|||_e^2\\[1ex]
& \qquad  + \alpha (\langle \omega(t_n), \nabla \cdot \vec E_{\tau,h}^1(t_n) \rangle - \langle \omega(t_{n-1}^+), \nabla \cdot \vec E_{\tau,h}^1(t_{n-1}^+) \rangle)\\[1ex]
& \qquad + \int_{I_n}\langle \vec T_{IV}^n, \Pi^{k-1}_\tau \vec E_{\tau,h}^2 \rangle \ud t  - \int_{I_n}  \langle \vec A_h \vec T_I^n, \Pi^{k-1}_\tau \vec E_{\tau,h}^2\rangle \ud t - \int_{I_n}  \langle  \vec T_{II}^n, \Pi^{k-1}_\tau \vec E_{\tau,h}^2 \rangle \ud t\\[1ex] 
& \qquad + \int_{I_n}  \langle \vec A_h \vec T^n_{III}, \Pi^{k-1}_\tau \vec E_{\tau,h}^2\rangle \ud t  - \alpha \int_{I_n}	\langle \partial_t \omega, \nabla \cdot \Pi^{k-1}_\tau \vec E_{\tau,h}^1 \rangle \ud t + \int_{I_n} \langle T_{V}^n, \Pi_\tau^{k-1} e_{\tau,h}\rangle \ud t \\[1ex]
& \qquad  - c_0 \int_{I_n} \langle \partial_t \omega, \Pi_\tau^{k-1} e_{\tau,h}\rangle \ud t  -\alpha  \int_{I_n} \langle \nabla \cdot \partial_t \vec \eta_1, \Pi_\tau^{k-1} e_{\tau,h}\rangle \ud t - \int_{I_n} \langle \vec K \nabla (p-I_\tau p), \nabla \Pi_\tau^{k-1} e_{\tau,h}\rangle \ud t\,.
\end{aligned}
\end{equation}
By the assumption of the positive-definiteness of $\vec K$, the inequalities of Cauchy--Schwarz and Cauchy--Young, identity \eqref{Eq:AhwAw} and integation by parts, applied to the last of the terms in \eqref{Eq:ES20}, we conclude from \eqref{Eq:ES20} that 
\begin{equation}
	\label{Eq:ES21}
	\begin{aligned}
		& ||| \vec E_{\tau,h} (t_n)|||_e^2 + |||e_{\tau,h}(t_n)|||_e^2 + c \int_{I_n} \|\nabla \Pi_\tau^{k-1} e_{\tau,h} \|^2 \ud t \leq ||| \vec E_{\tau,h} (t_{n-1}^+) |||_e^2 + |||e_{\tau,h}(t_{n-1}^+)|||_e^2  \\[1ex]
		& \qquad  + \delta_n - \delta_{n-1}^+  + c ||| \vec E_{\tau,h}|||^2_{L^2(I_n;\vec L^2)} + c \|e_{\tau,h}\|^2_{L^2(I_n;L^2)} + c \|\vec A \vec T_{I}^n\|_{L^2(I_n;\vec L^2)}^2  \\[1ex]
		& \qquad + c \left| \int_{I_n} \langle \vec T_{II}^n, \Pi^{k-1}_\tau \vec E_{\tau,h}^2 \rangle  \ud t \right|+c \|\vec A \vec T_{III}^n\|_{L^2(I_n;\vec L^2)}^2  + c \|\vec T_{IV}^n\|_{L^2(I_n;\vec L^2)}^2 +  c \| T_{V}^n\|_{L^2(I_n;L^2)}^2  \\[1ex] 
		& \qquad + c \left| \int_{I_n} \langle  \partial_t \omega, \Pi_\tau^{k-1} e_{\tau,h}\rangle \ud t\right| +  c \left| \int_{I_n} \langle \nabla \cdot \partial_t \vec \eta_1, \Pi_\tau^{k-1} e_{\tau,h}\rangle \ud t\right| + c \| p-I_\tau p\|_{L^2(I_n;H^2)}^2\,,
	\end{aligned}
\end{equation}
where $\delta_n$ and $\delta_{n-1}^+$ are defined in \eqref{Eq:ES00}. Combining \eqref{Eq:ES21} with Lem.~\ref{Lem:EstTni} and the bounds \eqref{Eq:ErrLI} and \eqref{Eq:Omega} proves the assertion \eqref{Eq:ES01} of this lemma.
\end{mproof}

Next, we estimate the right-hand side term $||| \vec E_{\tau,h} |||^2_{L^2(I_n;\vec L^2)} +\| e_{\tau,h}\|^2_{L^2(I_n;L^2)}$ in \eqref{Eq:ES01}.

\begin{lem}[Estimate of $||| \vec E_{\tau,h} |||^2_{L^2(I_n;\vec L^2)} +\| e_{\tau,h}\|^2_{L^2(I_n;L^2)}$]
\label{Lem:ErrL2L2}
Let $n=1,\ldots, N$. For the errors $\vec E_{\tau,h}$ and $e_{\tau,h}$, defined in \eqref{Eq:ErrU_00} and \eqref{Eq:ErrP_00}, there holds that 
\begin{equation}
\label{Eq:ErrL2L201}
\begin{aligned}
||| \vec E_{\tau,h} |||^2_{L^2(I_n;\vec L^2)} & +\| e_{\tau,h}\|^2_{L^2(I_n;L^2)}  \leq c  \tau_n (||| \vec E_{\tau,h}(t_{n-1}^+) |||^2 +\| e_{\tau,h}(t_{n-1}^+)\|^2)\\[1ex] 
& \quad + c \tau_n \big(\tau_n^{2(k+1)} (\mathcal E_t^{n,1})^2 + h^{2(r+1)} (\mathcal E_{\vec x}^{n,1}+\mathcal E_{\vec x}^{n,3})^2 + h^{2(r+2)} (\mathcal E_{\vec x}^{n,2})^2\big) \,,
\end{aligned} 
\end{equation}
where $\mathcal E_t^{n,1}$, $\mathcal E_{\vec x}^{n,1}$ and $\mathcal E_{\vec x}^{n,2}$ are defined in Lem.~\ref{Lem:ES} and $\mathcal E_{\vec x}^{n,3}:=\|p\|_{L^\infty(I_n;H^{r+1})}$.
\end{lem}

\begin{mproof}
Firstly, we consider \eqref{Eq:EE02}. Let $\vec E_{\tau,h}=(\vec E_{\tau,h}^1,\vec E_{\tau,h}^2)^\top$, defined in \eqref{Eq:ErrU_00}, be represented by 
\begin{equation}
\label{Eq:ErrL2L210}
\vec E_{\tau,h}^m(t) = \sum_{j=0}^k \vec E_{n,j}^{m} \phi_{n,j}(t)\,, \qquad \text{for}\; \;  t\in I_n\,, \; m\in\{1,2\}\,,
\end{equation}
where $\vec E_{n,j}^{m}\in \vec V_h^{r+1}$, for $j=0,\ldots,k$, and $\phi_{n,j}\in \mathbb P_k(I_n;\R)$, for $j=0,\ldots,k$, are the Lagrange interpolants with respect to $t_{n-1}$ and the Gauss quadrature nodes $t_{n,1}^{\text G},\ldots, t_{n,k}^{\text G}\in (t_{n-1},t_n)$ of \eqref{Eq:GF}. Then, it holds that $\vec E_{n,0}^{m}=\vec E_{\tau,h}^m(t_{n-1}^+)$. In \eqref{Eq:EE02}, we choose the test function 
\begin{equation}
\label{Eq:ErrL2L211}
	\vec \Phi_{\tau,h}(t) = \sum_{i=1}^k (\hat t_{i}^{\,\text G})^{-1/2} \begin{pmatrix}\vec A_h \vec{\widetilde E}_{n,i}^{1}\\[1ex] \vec{\widetilde E}_{n,i}^{2} \end{pmatrix} \psi_{n,i}(t)\,,
\end{equation}	
where $\vec{\widetilde E}_{n,i}^{m} := (\hat t_{i}^{\, \text G})^{-1/2} \vec{E}_{n,i}^{m}$, for $m\in \{1,2\}$ and $i=1,\ldots,k$, and $\psi_{n,i}\in \mathbb P_{k-1}(I_n;\R)$, for $i=1,\ldots,k$, are the Lagrange interpolants with respect to the Gauss quadrature nodes $t_{n,1}^{\text G},\ldots, t_{n,k}^{\text G}\in (t_{n-1},t_n)$ of \eqref{Eq:GF}. In \eqref{Eq:ErrL2L211}, the quantities $\hat t_{i}^{\, \text G}$, for $i=1,\ldots,k$, denote the quadrature nodes of the Gauss formula \eqref{Eq:GF} on the reference interval $\hat I$. Using the evaluation \eqref{Eq:ErrL2L210}, for the test function \eqref{Eq:ErrL2L211} it follows that 
\begin{equation}
\label{Eq:ErrL2L212}
\mbox{}\hspace*{-2cm}
\begin{aligned}
&\int_{I_n} \left\llangle \begin{pmatrix} \vec 0 & - \vec I\\ \vec A_h & \vec 0 \end{pmatrix} \begin{pmatrix} \vec{E}_{\tau,h}^1 \\ \vec{E}_{\tau,h}^2 \end{pmatrix}, \begin{pmatrix} \vec \Phi_{\tau,h}^1\\ \vec \Phi_{\tau,h}^2\end{pmatrix} \right\rrangle \ud t = \int_{I_n} \left\llangle \begin{pmatrix} \vec 0 & - \vec I\\ \vec A_h & \vec 0 \end{pmatrix} \begin{pmatrix} \vec{E}_{\tau,h}^1\\ \vec{E}_{\tau,h}^2 \end{pmatrix}, \sum_{i=1}^k (\hat t_{i}^{\,\text G})^{-1/2} \begin{pmatrix}\vec A_h \vec{\widetilde E}_{n,i}^{1}\\[1ex] \vec{\widetilde E}_{n,i}^{2} \end{pmatrix} \psi_{n,i}\right\rrangle \ud t\\[1ex]
& \qquad = \frac{\tau_n}{2} \sum_{\mu =1}^k \hat \omega_\mu^{\text G}\left\llangle \begin{pmatrix} -\vec{E}_{\tau,h}^2(t_{n,\mu}^{\text G})\\[1ex] \vec A_h \vec{E}_{\tau,h}^1(t_{n,\mu}^{\text G})  \end{pmatrix}, (\hat t_{\mu}^{\,\text G})^{-1/2} \begin{pmatrix}\vec A_h \vec{\widetilde E}_{n,\mu}^{1}\\[1ex] \vec{\widetilde E}_{n,\mu}^{2}\end{pmatrix}\right\rrangle\\[1ex]
& \qquad = \frac{\tau_n}{2} \sum_{\mu =1}^k \hat \omega_\mu^{\text G} (\hat t_{\mu}^{\,\text G})^{-1}\big(\langle -\vec{E}_{n,2}^\mu, \vec A_h \vec{E}_{\tau,h}^{1,\mu} \rangle + \langle \vec A_h \vec{E}_{n,1}^\mu, \vec{E}_{n,\mu}^{2}\rangle \big) = 0\,,
\end{aligned}\hspace*{-2cm}\mbox{}
\end{equation}
where the symmetry of $\vec A_h$ has been used in the last identity. By the expansion \eqref{Eq:ErrL2L211} along with the observation that $\vec E^m_{\tau,h}(t_{n-1}^+) = \vec E^m_{n,0}$, for $m\in \{1,2\}$, we have that 
\begin{equation}
	\label{Eq:ErrL2L213}	
	\begin{aligned}
		Q_n & := \int_{I_n} \left\llangle \begin{pmatrix} \vec I_d & \vec 0 \\ \vec 0 & \rho\vec I_d \end{pmatrix}\begin{pmatrix} \partial_t \vec{E}_{\tau,h}^1\\ \partial_t \vec{E}_{\tau,h}^2 \end{pmatrix},\begin{pmatrix} \vec \Phi_{\tau,h}^1\\ \vec \Phi_{\tau,h}^2\end{pmatrix} \right\rrangle \ud t\\[1ex] 
		& = \int_{I_n} \left\llangle \begin{pmatrix} 1& 0 \\ 0 & \rho \end{pmatrix}\begin{pmatrix} \partial_t \vec{E}_{\tau,h}^1\\ \partial_t \vec{E}_{\tau,h}^2 \end{pmatrix}, \sum_{i=1}^k (\hat t_{i}^{\,\text G})^{-1/2}  \begin{pmatrix}\vec A_h \vec{\widetilde E}_{n,i}^{1}\\[1ex] \vec{\widetilde E}_{n,i}^{2} \end{pmatrix} \psi_{n,i}\right\rrangle \ud t\\[1ex]
		 & = \sum_{i,j=1}^k \widetilde m_{ij} \big(\langle \vec A_h \vec{\widetilde E}_{n,j}^{1}, \vec{\widetilde E}_{n,i}^{1}\rangle + \rho \langle \vec{\widetilde E}_{n,j}^{2}, \vec{\widetilde E}_{n,i}^{2}\rangle\big)\\[1ex] 
		 &  \qquad + \sum_{i=1}^k m_{i0} (\hat t_{i}^{\,\text G})^{-1/2}  \big(\langle \vec A_h \vec{E}_{\tau,h}^{1}(t_{n-1}^+), \vec{\widetilde E}_{n,i}^{1}\rangle + \rho \langle \vec{E}_{\tau,h}^{2}(t_{n-1}^+), \vec{\widetilde E}_{n,i}^{2}\rangle\big)\,,
	\end{aligned}
\end{equation}
where the matrix $\vec{M} = (m_{ij})_{i,j=1,\ldots,k}$ and vector $\vec m_0=(m_{i0})_{i=1,\ldots,k}$ are defined by 
\begin{equation*}
	\begin{aligned}
	m_{ij}  := \int_{I_n} \phi^{\prime}_{n,j}(t) \psi_{n,i}(t) \ud t \,, \;\; \text{for} \;\; i\in\{1,\ldots, k\}\,, \; j \in \{1,\ldots, k\}\,,\quad 
	m_{i0}  := \int_{I_n} \phi^{\prime}_{n,0}(t) \psi_{n,i}(t) \ud t \,, \;\; \text{for} \;\; i\in\{1,\ldots, k\}\,,
	\end{aligned}
\end{equation*}
and the matrix $\vec{\widetilde M} = (\widetilde m_{ij})_{i,j=1,\ldots,k}$ is given by 
\begin{equation*}
\vec{\widetilde M}:=\vec D^{-1/2}\vec M \vec D^{1/2}\,, \quad \text{with} \;\; \vec D = \operatorname{diag}\{\hat t_{1}^{\text G}, \ldots,\hat t_{k}^{\text G}\}\,.
\end{equation*}
By the positivity of $\vec{\widetilde M}$ (cf.~\cite[Lem.~2.1]{KM04}) we then have that 
\begin{equation} 
\label{Eq:ErrL2L214}	
Q_n \geq c \sum_{j=1}^k ||| \vec{\widetilde E}_{n,j} |||^2 - c \Big(\sum_{j=1}^k ||| \vec{\widetilde E}_{n,j} |||^2 \Big)^{1/2}||| \vec{E}_{\tau,h}(t_{n-1}^+) ||| \geq c \sum_{j=1}^k ||| \vec{\widetilde E}_{n,j} |||^2 - c ||| \vec{E}_{\tau,h}(t_{n-1}^+) |||^2\,.
\end{equation}
By the equivalence of $\sum_{j=1}^k ||| \vec{\widetilde E}_{n,j} |||$ and $\sum_{j=1}^k ||| \vec{E}_{n,j}|||$ along with the equivalence (cf.~\cite[Eq.~(2.4)]{KM04})
\begin{equation} 
\label{Eq:ErrL2L215}	
c_1 \tau_n \sum_{j=0}^k ||| \vec{E}_{n,j}|||^2 \leq ||| \vec E_{\tau,h} |||^2_{L^2(I_n;\vec L^2)}\leq c_2 \tau_n \sum_{j=0}^k ||| \vec{E}_{n,j}|||^2 \,, 
\end{equation} 
we conclude from \eqref{Eq:ErrL2L213} to \eqref{Eq:ErrL2L215} that 
\begin{equation} 
\label{Eq:ErrL2L216}	
\tau_n Q_n \geq c ||| \vec{E}_{\tau,h}|||^2_{L^2(I_n;\vec L^2)} -  c \tau_n ||| \vec{E}_{\tau,h}(t_{n-1}^+)|||^2\,.
\end{equation} 
Next, we address the last term on the left-hand side of \eqref{Eq:EE02} for the test function \eqref{Eq:ErrL2L211}. Similarly to \eqref{Eq:ErrL2L210}, for the error $e_{\tau,h}$ we use the representation 
\begin{equation}
	\label{Eq:ErrL2L250}
	e_{\tau,h}(t) = \sum_{j=0}^k e_{n,j} \phi_{n,j}(t)\,, \qquad \text{for}\; \;  t\in I_n\,,
\end{equation}
where $e_{n,j}\in V_h^r$, for $j=0,\ldots,k$. Further we put $\widetilde{e}_{n,i}:=(\hat t_i^{\text G})^{-1/2}{e}_{n,i}$, for $i=1,\ldots,k$. Using \eqref{Eq:ErrL2L211} and \eqref{Eq:ErrL2L250} along with \eqref{Eq:ES13} and recalling that $\psi_{n,i}\in \mathbb P_{k-1}(I_n;\R)$ in \eqref{Eq:ErrL2L211}, it follows that 
\begin{equation}
\label{Eq:ErrL2L217}	
	\begin{aligned}
		\int_{I_n} \left\langle e_{\tau,h}, \nabla \cdot \vec \Phi_{\tau,h}^2 \right\rangle \ud t & = \int_{I_n} \left\langle e_{\tau,h} , \nabla \cdot \sum_{i=1}^k (\hat t_{i}^{\,\text G})^{-1/2}  \vec{\widetilde E}_{n,i}^{2}\psi_{n,i}\right\rangle \ud t \\[1ex] 
		& = \frac{\tau_n}{2} \sum_{\mu=1}^k \hat \omega_\mu^{\text G}  \left\langle e_{\tau,h}(t_{n,\mu}^{\text G}) , \nabla \cdot \sum_{i=1}^k (\hat t_{i}^{\,\text G})^{-1/2}  \vec{\widetilde E}_{n,i}^{2}\psi_{n,i}(t_{n,\mu}^{\text G})\right\rangle\\[1ex]
		& =  \frac{\tau_n}{2} \sum_{\mu=1}^k \hat \omega_\mu^{\text G} (\hat t_{\mu}^{\,\text G})^{-1}  \langle e_{n,\mu} , \nabla \cdot \vec{E}_{\tau,h}^2(t_{n,\mu}^{\text G})\rangle\\[1ex] 
		& =  \frac{\tau_n}{2} \sum_{\mu=1}^k \hat \omega_\mu^{\text G} (\hat t_{\mu}^{\,\text G})^{-1}  \langle e_{n,\mu} , \nabla \cdot \partial_t \vec{E}_{\tau,h}^1(t_{n,\mu}^{\text G})\rangle\\[1ex] 
		& = \int_{I_n} \Big\langle \nabla \cdot \partial_t \vec{E}_{\tau,h}^1,\sum_{i=1}^k (\hat t_i^{\text G})^{-1/2} \widetilde{e}_{n,i} \psi_{n,i}\Big\rangle \ud t \,.
\end{aligned}
\end{equation}
Finally, we address the last term on the right-hand side of \eqref{Eq:EE02} with \eqref{Eq:ErrL2L211}. Similarly to \eqref{Eq:ErrL2L217}, using \eqref{Eq:ErrL2L211} with $\psi_{n,i}\in \mathbb P_{k-1}(I_n;\R)$ and employing \eqref{Eq:ES13}, we find that 
\begin{equation}
\label{Eq:ErrL2L218}	
\begin{aligned}
R_n & := \int_{I_n} \Big\langle \omega, \nabla \cdot \vec \Phi_{\tau,h}^2 \Big\rangle \ud t  = \int_{I_n} \Big\langle \Pi_\tau^{k-1}\omega , \nabla \cdot \sum_{i=1}^k (\hat t_{i}^{\,\text G})^{-1/2}  \vec{\widetilde E}_{n,i}^{2}\psi_{n,i}\Big\rangle \ud t\\[1ex]
&= \frac{\tau_n}{2} \sum_{\mu=1}^k \hat \omega_\mu^{\text G} (\hat t_{\mu}^{\,\text G})^{-1}  \langle  \Pi_\tau^{k-1} \omega(t_{n,\mu}^{\text G}), \nabla \cdot \vec{E}_{\tau,h}^2(t_{n,\mu}^{\text G}) \rangle\\[1ex]
& = \frac{\tau_n}{2} \sum_{\mu=1}^k \hat \omega_\mu^{\text G}  \Big\langle \sum_{i=1}^k (\hat t_{i}^{\,\text G})^{-1} \Pi_\tau^{k-1}\omega(t_{n,i}^{\text G}) \psi_{n,i}(t_{n,\mu}^{\text G}), \nabla \cdot \partial_t  \vec{E}_{\tau,h}^1(t_{n,\mu}^{\text G}) \Big\rangle\\[1ex]
& = \int_{I_n} \Big\langle \sum_{i=1}^k (\hat t_{i}^{\,\text G})^{-1} \Pi_\tau^{k-1}\omega(t_{n,i}^{\text G}) \psi_{n,i}, \nabla \cdot \partial_t  \vec{E}_{\tau,h}^1\Big\rangle  \ud t \,.
\end{aligned}
\end{equation}
From \eqref{Eq:ErrL2L218} along with $\int_{I_n} \psi_{n,i}^2 \ud t \leq c \tau_n$ and the inequality of Cauchy--Young we get that 
\begin{equation*}
\label{Eq:ErrL2L219}	
\begin{aligned}
\tau_n R_n & \leq c \tau_n \max_{i=1,\ldots,k} \{\| \Pi_\tau^{k-1}\omega(t_{n,i}^{\text G})\| \, \|\psi_{n,i}\|_{L^2(I_n;\R)}\} \| \partial_t \nabla \cdot \vec{E}_{\tau,h}^1 \|_{L^2(I_n;L^2)} \\[1ex]
&\leq c \tau_n \max_{i=1,\ldots,k} \| \Pi_\tau^{k-1}\omega(t_{n,i}^{\text G})\|^2 + \varepsilon \tau_n^2 \| \partial_t \nabla \cdot \vec{E}_{\tau,h}^1 \|_{L^2(I_n; L^2)}^2
\end{aligned}
\end{equation*}
with a sufficiently small constant $\varepsilon >0$. The $L^\infty$--$L^2$ inverse relation \eqref{Eq:LinfL2}, the error estimate \eqref{Eq:Omega} for the elliptic projection $R_h$ in $\omega=p-R_h p$ and the $H^1$--$L^2$ inverse inequality then imply that 
\begin{equation}
\label{Eq:ErrL2L220}	
\begin{aligned}
\tau_n R_n & \leq c \|\Pi_\tau^{k-1} w\|^2_{L^2(I_n;L^2)} + \varepsilon ||| \vec{E}_{\tau,h} ||| _{L^2(I_n; L^2)}^2 \leq c \tau_n \| w\|^2_{L^\infty(I_n;L^2)} + \varepsilon ||| \vec{E}_{\tau,h} ||| _{L^2(I_n; L^2)}^2\\[1ex]
& \leq c  \tau_n h^{2(r+1)} \| p \|^2_{L^\infty(I_n;H^{r+1})} + \varepsilon ||| \vec{E}_{\tau,h} ||| _{L^2(I_n; L^2)}^2\,.
\end{aligned}
\end{equation}
For a suitable choice of $\varepsilon$, the second term on right-hand side of \eqref{Eq:ErrL2L220} can be absorbed by the left-hand side of \eqref{Eq:ErrL2L201}. The remaining terms on the right-hand side of \eqref{Eq:EE02} can be treated as before in Lem.~\ref{Lem:ES}. 

Now, we consider \eqref{Eq:EE03}. We choose the test function 
\begin{equation}
\label{Eq:ErrL2L251}
\psi_{\tau,h}(t)  =  \sum_{i=1}^k (\hat t_i^{\text G})^{-1/2} \widetilde{e}_{n,i} \psi_{n,i}(t)\,.	
\end{equation}
By arguments similarly to \eqref{Eq:ErrL2L213} to \eqref{Eq:ErrL2L216} and with \eqref{Eq:ErrL2L250}, we then have that 
\begin{equation}
\label{Eq:ErrL2L252}
\begin{aligned}
S_n & := \int_{I_n} \langle \partial_t e_{\tau,h}, \psi_{\tau,h}\rangle \ud t  =  \int_{I_n} \Big\langle \partial_t e_{\tau,h}, \sum_{i=1}^k (\hat t_i^{\text G})^{-1/2} \widetilde{e}_{n,i} \psi_{n,i}\Big \rangle \ud t\\[1ex]
&  = \dfrac{\tau_n}{2} \sum_{\mu=1}^k \hat \omega_\mu^{\text G}\Big\langle \partial_t e_{\tau,h}(t_{n,\mu}^{\text G}), \sum_{i=1}^k (\hat t_i^{\text G})^{-1/2} \widetilde{e}_{n,i} \psi_{n,i}(t_{n,\mu}^{\text G})\Big\rangle \ud t \\[1ex]
& = \sum_{i,j=1}^k \widetilde m_{ij} \langle \widetilde e_{n,j}, \widetilde e_{n,i}\rangle + \sum_{i=1}^k m_{i0} (\hat t_{i}^{\,\text G})^{-1/2} \langle {e}_{\tau,h}(t_{n-1}^+), \widetilde e_{n,i}\rangle\\[1ex]
&  \geq c \sum_{j=1}^k \| \widetilde e_{n,j} \|^2_{L^2(I_n;L^2)} - c \| e_{\tau,h}(t_{n-1}^+\|^2\,.
\end{aligned}
\end{equation}
Similarly to \eqref{Eq:ErrL2L216}, we conclude from \eqref{Eq:ErrL2L252} that
\begin{equation}
\label{Eq:ErrL2L2525}
\tau_n S_n \geq c \|e_{\tau,h}\|^2_{L^2(I_n;L^2)} - c \tau_n \| e_{\tau,h}(t_{n-1}^+\|^2\,.
\end{equation}
Further, we obtain by \eqref{Eq:ErrL2L250} along with the positive definiteness of $\vec K$ that 
\begin{equation}
	\label{Eq:ErrL2L253}
	\begin{aligned}
		\int_{I_n} \langle B_h e_{\tau,h}, \psi_{\tau,h}\rangle \ud t & = \int_{I_n} \Big\langle B_h e_{\tau,h}, \sum_{i=1}^k (\hat t_i^{\text G})^{-1/2} \widetilde{e}_{n,i} \psi_{n,i} \Big\rangle \ud t \\[1ex]
		&  = \dfrac{\tau_n}{2} \sum_{\mu=1}^k \hat \omega_\mu^{\text G} \Big\langle B_h e_{\tau,h}(t_{n,\mu}^{\text G}), \sum_{i=1}^k (\hat t_i^{\text G})^{-1/2} \widetilde{e}_{n,i} \psi_{n,i}(t_{n,\mu}^{\text G}) \Big\rangle\\[1ex]
		& = \dfrac{\tau_n}{2} \sum_{\mu=1}^k \hat \omega_\mu^{\text G} (\hat t_i^{\text G})^{-1}  \Big\langle B_h e_{\tau,h}(t_{n,\mu}^{\text G}), e_{\tau,h}(t_{n,\mu}^{\text G}) \Big\rangle \geq 0\,.
	\end{aligned}
\end{equation}
The terms on the right-hand side of \eqref{Eq:EE03} can be treated as before in Lem.~\ref{Lem:ES}. 

Finally, we sum up the error equations \eqref{Eq:EE02} and \eqref{Eq:EE03} for the test functions \eqref{Eq:ErrL2L211} and \eqref{Eq:ErrL2L251}. After summation, we use \eqref{Eq:ErrL2L212}, \eqref{Eq:ErrL2L213}, \eqref{Eq:ErrL2L216}, \eqref{Eq:ErrL2L217}, \eqref{Eq:ErrL2L218} and \eqref{Eq:ErrL2L220} along with \eqref{Eq:ErrL2L252}, \eqref{Eq:ErrL2L2525} and \eqref{Eq:ErrL2L253}. The remaining terms are treated as before in the proof of Lem.~\ref{Lem:ES}. By \eqref{Eq:ErrL2L217}, the terms $-\alpha\int_{I_n} \left\langle e_{\tau,h}, \nabla \cdot \vec \Phi_{\tau,h}^2 \right\rangle \ud t$ and $\alpha\int_{I_n} \Big\langle \nabla \cdot \partial_t \vec{E}_{\tau,h}^1,\psi_{\tau,h}\rangle \ud t$ cancel out for the test functions \eqref{Eq:ErrL2L211} and \eqref{Eq:ErrL2L251}. This is one of the key arguments of our proof for controling the coupling mechanism of the subsystems. By using Lem.~\ref{Lem:EstTni} and employing the inequality of Cauchy--Young, we then conclude the assertion \eqref{Eq:ErrL2L201} of this lemma. 
\end{mproof}

It remains to estimate $||| \vec E_{\tau,h}(t_{n-1}^+) |||_e^2 + ||| e_{\tau,h}(t_{n-1}^+)|||_e^2$, arising on the right-hand side of \eqref{Eq:ES01}. 

\begin{lem}[Estimate of $||| \vec E_{\tau,h}(t_{n-1}^+) |||_e^2 +|||e_{\tau,h}(t_{n-1}^+)|||_e^2$]
\label{Lem:En_1+}
Let $n=\in\{2,\ldots, N\}$. For the errors $\vec E_{\tau,h}$ and $e_{\tau,h}$, defined in \eqref{Eq:ErrU_00} and \eqref{Eq:ErrP_00}, there holds that 
\begin{equation}
\label{Eq:En_1+01}
\begin{aligned}
||| \vec E_{\tau,h}(t_{n-1}^+) |||_e^2 + ||| e_{\tau,h}(t_{n-1}^+)|||_e^2 & \leq (1+\tau_{n-1}) (||| \vec E_{\tau,h}(t_{n-1}) |||_e^2 + ||| e_{\tau,h}(t_{n-1})|||_e^2)  + c \tau_{n-1}^{2(k+1)} (\mathcal E_{t}^{n-1,2})^2 \,,
\end{aligned}
\end{equation}
where $\mathcal  E_{t}^{n-1,2}:= \| \partial_t^{k+2} \vec u\|_{L^2(I_{n-1};\vec H^1)}$.
\end{lem}

\begin{mproof}
Firstly, from \eqref{Eq:ErrU_00}, \eqref{Eq:DefW}, the continuity constraint imposed on $\vec v_{\tau,h}$ in Problem \ref{Prob:DPV} or \ref{Prob:DPQ}, respectively, and the assumption that $\vec u$ is sufficient regular we get that 
\begin{equation}
\label{Eq:En_1+10}
\begin{aligned}
\vec E_{\tau,h}^2(t_{n-1}^+) & = \vec w_2(t_{n-1}^+) - \vec v_{\tau,h}(t_{n-1}^+)=  \vec R_h \partial_t \vec u (t_{n-1}^+) - \vec v_{\tau,h} (t_{n-1}^+) \\[1ex]
& =  \vec R_h \partial_t \vec u (t_{n-1}) - \vec v_{\tau,h} (t_{n-1}) = \vec E_{\tau,h}^2(t_{n-1})\,.
\end{aligned}
\end{equation}
For \eqref{Eq:En_1+10}, we recall the notation that $\vec v_{\tau,h} (t_{n-1}^+) = \lim_{t\rightarrow t_{n-1}}\vec v_{\tau,h}{}_{|I_n} (t)$ and $\vec v_{\tau,h} (t_{n-1}) = \vec v_{\tau,h}{}_{|I_{n-1}} (t_{n-1})$. Secondly, by the continuity constraint imposed on $\vec u_{\tau,h}$ in Problem \ref{Prob:DPV} or \ref{Prob:DPQ}, respectively, we have that 
\begin{equation*}
\vec E_{\tau,h}^1(t_{n-1}^+) =  \vec w_1(t_{n-1}^+) - \vec u_{\tau,h}(t_{n-1}^+) = (\vec w_1(t_{n-1}^+)  - \vec w_1(t_{n-1})) +  \vec E^1_{\tau,h}(t_{n-1}) \,.
\end{equation*}
Then, by the triangle inequality of the norm property of \eqref{Eq:DefTNe} along with \eqref{Eq:Nequi} we can conclude that
\begin{equation}
\label{Eq:En_1+11}
\begin{aligned}
\langle \vec C & \vec \varepsilon(\vec E_{\tau,h}^1(t_{n-1}^+)), \vec \varepsilon(\vec E_{\tau,h}^1(t_{n-1}^+)) \rangle^{1/2} \\[1ex]
& \leq c \| \nabla (\vec w_1(t_{n-1}^+)  - \vec w_1(t_{n-1}))\| + \langle \vec C \vec \varepsilon(\vec E_{\tau,h}^1(t_{n-1})), \vec \varepsilon(\vec E_{\tau,h}^1(t_{n-1}))\rangle^{1/2} \,.
\end{aligned}
\end{equation}
By \eqref{Eq:DefW} and \eqref{Eq:ErrLI} there holds that 
\begin{equation}
\label{Eq:En_1+12}
\begin{aligned}
\| \nabla (\vec w_1(t_{n-1}^+)  - \vec w_1(t_{n-1}))\|  & \leq c \| \nabla (\vec w_1(t_{n-1}^+)  - \vec R_h  \vec u(t_{n-1}))\|\\[1ex]
& \qquad +  \Big\| \nabla\vec R_h  \Big(\vec u(t_{n-1}) - \int_{t_{n-2}}^{t_{n-1}} I_\tau (\partial_t u) \ud t - \vec u({t_{n-2}})\Big)\Big\| \\[1ex]
& \leq c  \Big\| \int_{t_{n-2}}^{t_{n-1}} \nabla \partial_t \vec u - I_\tau (\nabla \partial_t \vec u) \ud t\Big\| \leq c \tau_{n-1}^{1/2} \, \tau_{n-1}^{k+1} \, \| \partial_t^{k+2} \vec u\|_{L^2(I_{n-1};\vec H^1)}\,.
\end{aligned}
\end{equation}
Thirdly, since $ e_{\tau,h} \in X_\tau^k(V_h^r) \subset C([0,T];V_h^r)$ we have that 
\begin{equation}
\label{Eq:En_1+20}
e_{\tau,h}(t_{n-1}^+)  =  e_{\tau,h}(t_{n-1}) \,.
\end{equation}
Combining \eqref{Eq:En_1+10} to \eqref{Eq:En_1+20} and applying the arithmetic and geometric mean inequality proves the assertion \eqref{Eq:En_1+01}. 
\end{mproof}

The term $||| \vec E_{\tau,h}(t_{n-1}^+) |||^2 + \| e_{\tau,h}(t_{n-1}^+)\|^2$,  arising on the right-hand side of \eqref{Eq:ErrL2L201}, can be estimated along the lines of \eqref{Eq:En_1+01} as well. Finally, we address the term $\delta_n-\delta_{n-1}^+$ of \eqref{Eq:ES01}. 

\begin{lem}
\label{Lem:DifDelta}
Let $\delta_n$ and $\delta_{n-1}^+$ be defined by \eqref{Eq:ES00}. For $n=2,\ldots,N$ there holds that
\begin{equation}
\label{Eq:DD_01}
\begin{aligned}
\delta_n-\delta_{n-1}^+ \leq  \delta_n-\delta_{n-1}  + 	c  \tau_{n-1} \, \tau_{n-1}^{2(k+1)} \, (\mathcal E_{t}^{n-1,3})^2 + c \tau_{n-1}\, h^{2(r+1)}\, (\mathcal E_{\vec x}^{n-1,3})^2 \,,
\end{aligned}
\end{equation}
where $\mathcal E_{t}^{n-1,3} := \| \partial_t^{k+2} \vec u \|_{L^\infty(I_{n-1};\vec H^1)}$ and $\mathcal E_{\vec x}^{n-1,3} $ is defined by Lem.~\ref{Lem:ErrL2L2}. For $n=1$ there holds that
\begin{equation}
	\label{Eq:DD_02}
	\begin{aligned}
		|\delta_1-\delta_{0}^+|  \leq c h^{2(r+1)}  \big(\|p_0\|^2_{r+1}+ \|p(t_1)\|^2_{r+1}+\| \vec u_0\|^2_{r+2}\big) + \varepsilon ||| \vec E_{\tau,h}(t_1)|||_e^2
	\end{aligned}
\end{equation}
for a (sufficiently small) constant $\varepsilon >0$.
\end{lem}

\begin{mproof}
By definition \eqref{Eq:ES00} of $\delta_n$ and $\delta_{n-1}^+$ along with \eqref{Eq:ErrU_00}, \eqref{Eq:ErrP_00},  \eqref{Eq:DefW}, the interpolation property \eqref{Eq:DefLagIntOp} of $I_\tau$, the continuity of $\vec u_{\tau,h}$ and the approximation properties \eqref{Eq:ErrLI} and \eqref{Eq:Omega} we have for $n=2,\ldots,N$ that 
\begin{equation*}
	\begin{aligned}
		\delta_{n-1}^+ & = \alpha \langle p (t_{n-1})- R_h p(t_{n-1}), \nabla \cdot (\vec R_h \vec u(t_{n-1}) - \vec u_{\tau,h}(t_{n-1}))\rangle\\[1ex]
		& =  \alpha \langle p (t_{n-1})- R_h p(t_{n-1}), \nabla \cdot (\vec w_1(t_{n-1}) - \vec u_{\tau,h}(t_{n-1}))\rangle\\[1ex]
		&\qquad  + \alpha \langle p (t_{n-1})- R_h p(t_{n-1}), \nabla \cdot (\vec R_h \vec u(t_{n-1}) - \vec w_1(t_{n-1}))\rangle\\[1ex]
		& = \delta_{n-1} + \alpha \Big\langle p (t_{n-1})- R_h p(t_{n-1}), \nabla \cdot \Big(\vec R_h \vec u(t_{n-1}) -\int_{t_{n-2}}^{t_{n-1}} I_\tau (\vec R_h\partial_t \vec u) \ud t - \vec R_h \vec u(t_{n-2}) \Big)\Big\rangle \\[1ex]
		& = \delta_{n-1} + \alpha \Big\langle p (t_{n-1})- R_h p(t_{n-1}), \nabla \cdot \int_{t_{n-2}}^{t_{n-1}} (\vec R_h\partial_t \vec u) - I_\tau (\vec R_h\partial_t \vec u) \ud t \big\rangle\\[1ex] 
		&  = \delta_{n-1} + \alpha \Big\langle p (t_{n-1})- R_h p(t_{n-1}), \nabla \cdot \int_{t_{n-2}}^{t_{n-1}} \vec R_h(\partial_t \vec u - I_\tau (\partial_t \vec u)) \ud t \big\rangle =: \delta_{n-1} + \varepsilon_{n-1}\,,
	\end{aligned}
\end{equation*}
such that 
\begin{equation}
	\label{Eq:DD_04}
\delta_n-\delta_{n-1}^+ = \delta_n-\delta_{n-1} - \varepsilon_{n-1}\,,
\end{equation}
where 
\begin{equation}
	\label{Eq:DD_03}
	\begin{aligned}
		|\varepsilon_{n-1}| & \leq c \tau_{n-1}\, h^{r+1}\,\mathcal E_{\vec x}^{n-1,3} \, \tau_{n-1}^{k+1} \, \mathcal E_{t}^{n-1,3}  \leq c \tau_{n-1} \, \tau_{n-1}^{2(k+1)} \, (\mathcal E_{t}^{n-1,3})^2 +c \tau_{n-1}\, h^{2(r+1)}\, (\mathcal E_{\vec x}^{n-1,3})^2  \,. 
	\end{aligned}
\end{equation}
Now, the assertion \eqref{Eq:DD_01} is a direct consequence of \eqref{Eq:DD_04} and \eqref{Eq:DD_03}.

For $n=1$, there holds by \eqref{Eq:ES00}, \eqref{Eq:ErrU_00}, \eqref{Eq:ErrP_00} and \eqref{Eq:DefW} along with the Assumption \ref{AssIV} that 
\begin{equation}
	\label{Eq:DD_05}
	\begin{aligned}
		\delta_{0}^+ = \alpha \langle p_0 - R_h p_0, \nabla \cdot (\vec R_h \vec u_0 - \vec u_{0,h})\rangle \leq c h^{2(r+1)}\big(\|p_0\|^2_{r+1}+ \| \vec u_0\|^2_{r+2}\big)\,.
	\end{aligned}
\end{equation}
Further, by the inequalities of Cauchy--Schwarz and Cauchy--Young along with \eqref{Eq:Omega} we have that 
\begin{equation}
	\label{Eq:DD_06}	
	\begin{aligned}
		\delta_1 & = \alpha \langle p (t_{1})- R_h p(t_{1}), \nabla \cdot \vec E^1_{\tau,h}(t_1)\rangle \leq c h^{2(r+1)}\|p(t_1)\|^2_{r+1} + \tilde \varepsilon \, ||| \vec E_{\tau,h}(t_1) |||^2
	\end{aligned}
\end{equation}
with $\tilde \varepsilon>0$. By \eqref{Eq:Nequi} and the triangle inequality we get \eqref{Eq:DD_02} from \eqref{Eq:DD_05} and \eqref{Eq:DD_06}. 
\end{mproof}

\begin{thm}[Main convergence result]
\label{Thm:MainRes}	
For the approximation $(\vec u_{\tau,h},\vec v_{\tau,h},p_{\tau,h})$ defined by Problem \eqref{Prob:DPV} or \eqref{Prob:DPQ}, respectively, of the sufficiently regular solution $(\vec u,\vec v,  p)$ with $\vec v = \partial_t \vec u$ to \eqref{Eq:HPS} there holds that
 \begin{equation}
\label{Eq:MR01}
\|\nabla (\vec u(t) - \vec u_{\tau,h} (t))\|  + \| \vec v(t) - \vec v_{\tau,h} (t)\|  + \|p(t) - p_{\tau,h}(t)\| \leq c \tau^{k+1}+ c h^{r+1}\,, \quad \text{for} \;\; t\in I\,. 
\end{equation}

\end{thm}

\begin{mproof}
Combining the estimates \eqref{Eq:ES01} and \eqref{Eq:ErrL2L201} and recalling the norm equivalence \eqref{Eq:Nequi} yields that 
\begin{equation}
\label{Eq:MR10}
\begin{aligned}
 ||| \vec E_{\tau,h} (t_n)|||_e^2  & + |||e_{\tau,h}(t_n)|||_e^2   \leq \delta_n - \delta_{n-1}^+ + (1+c\tau_n ) (||| \vec E_{\tau,h} (t_{n-1}^+) |||_e^2 + |||e_{\tau,h}(t_{n-1}^+)|||_e^2) \\[1ex]
& \qquad + c \big(\tau_n^{2(k+1)} (\mathcal E_t^{n,1})^2 + ch^{2(r+1)} (\mathcal E_{\vec x}^{n,1})^2 + h^{2(r+2)} (\mathcal E_{\vec x}^{n,2})^2\big)  + c\tau_n h^{2(r+1)} (\mathcal E_{\vec x}^{n,3})^2 
\end{aligned}
\end{equation}
for  $n=1,\ldots,N$. Employing now \eqref{Eq:En_1+01} and \eqref{Eq:DD_01} in \eqref{Eq:MR10}, implies that
\begin{equation}
	\label{Eq:MR11}
	\begin{aligned}
		 ||| \vec E_{\tau,h} (t_n)|||_e^2  & + |||e_{\tau,h}(t_n)|||_e^2   \leq \delta_n - \delta_{n-1}+ (1+c\tau_n)(1+\tau_{n-1}) \big(||| \vec E_{\tau,h}(t_{n-1}) |||_e^2 + |||e_{\tau,h}(t_{n-1})|||_e^2\big)\\[1ex] 
		& + c \big(\tau_n^{2(k+1)} (\mathcal E_t^{n,1})^2 + h^{2(r+1)} (\mathcal E_{\vec x}^{n,1})^2 + c h^{2(r+2)} (\mathcal E_{\vec x}^{n,2})^2\big) + \tau_n h^{2(r+1)} (\mathcal E_{\vec x}^{n,3})^2	  \\[1ex] 
		& + c \tau_{n-1}^{2(k+1)} (\mathcal E_{t}^{n-1,2})^2 + \tau_{n-1}\tau_{n-1}^{2(k+1)}(\mathcal E_{t}^{n-1,3})^2 + \tau_{n-1}h^{2(r+1)}(\mathcal E_{\vec x}^{n-1,3})^2
	\end{aligned}
\end{equation}
for $n=2,\ldots, N$.  It remains to consider the case that $n=1$. By Problem~\ref{Prob:DPV} we have that $\vec U_{\tau,h}\in (C([0,T];$ $V_h^{r+1}))^{2d}$ and $p_{\tau,h}\in C([0,T];V_h^r)$. By \eqref{Eq:DefW} we have that $\vec w_1(t_0)= \vec R_h \vec u_0$ and $\vec w_2(t_0)= \vec R_h \vec u_1$. Thus, for $||| \vec E_{\tau,h} (t_0^+)|||_e$ and $ ||| e_{\tau,h}(t_0^+)|||_e$ it follows under the Assumption~\ref{AssIV} that 
\begin{equation}
\label{Eq:MR12}
\begin{aligned}
||| \vec E_{\tau,h} (t_0^+)|||_e^2  + ||| e_{\tau,h}(t_0^+)|||^2  & \leq c  \| \nabla (\vec R_h \vec u_0 - \vec u_{0,h})\|^2 + c \| \vec R_h \vec u_1 - \vec v_{0,h}\|^2 +c \| R_h p_0 - p_{0,h}\|^2 \leq c  h^{2(r+1)}\,. 
\end{aligned}
\end{equation}
Employing \eqref{Eq:MR12} and \eqref{Eq:DD_02} in \eqref{Eq:MR10}, we obtain that, for sufficiently regular solutions $(\vec u,p)$ \eqref{Eq:HPS},  
\begin{equation}
	\label{Eq:MR13}
	\begin{aligned}
		& ||| \vec E_{\tau,h} (t_1)|||_e^2  + |||e_{\tau,h}(t_1)|||_e^2   \leq c \tau_1^{2(k+1)} + c h^{2(r+1)}\,.
	\end{aligned}
\end{equation}

Next, we introduce the abbreviation that
\begin{equation}
\label{Eq:MR20}	
A_n :=  ||| \vec E_{\tau,h} (t_n)|||_e^2  + |||e_{\tau,h}(t_n)|||_e^2\,, \quad \text{for}\;\; n=0,\ldots, N\,.
\end{equation}
Then, we recover \eqref{Eq:MR11} as 
\begin{equation}
\label{Eq:MR21}	
	\begin{aligned}
	A_n & \leq \delta_n - \delta_{n-1} + (1+c\tau_n)(1+\tau_{n-1}) A_{n-1} + c\tau_n^{2(k+1)} \big((\mathcal E_t^{n,1})^2 + (\mathcal E_t^{n-1,2})^2 \big)\\[1ex] 
	& \qquad  + c h^{2(r+1)} \big((\mathcal E_{\vec x}^{n,1})^2 + (\mathcal E_{\vec x}^{n,2})^2\big) + c\tau_n\, \tau_n^{2(k+1)} (\mathcal E_{t}^{n-1,3})^2	 + c \tau_n\, h^{2(r+1)}\, \big((\mathcal E_{\vec x}^{n,3})^2+(\mathcal E_{\vec x}^{n-1,3})^2 \big) 
\end{aligned}
\end{equation}
for $n=2,\ldots,N$. From \eqref{Eq:MR13} we have that
\begin{equation}
	\label{Eq:MR23}	
	\begin{aligned}
		A_1 & \leq c_1 \tau_1^{2(k+1)}+ c_2 h^{2(r+1)}\,.
\end{aligned}
\end{equation}
Now, we apply the discrete Gronwall inequality \cite[Lem.~1.4.2]{QV94} to \eqref{Eq:MR21} and \eqref{Eq:MR23}.  For this,  we change the index $n$ to $m$ in \eqref{Eq:MR21} and sum up the resulting inequality from $m=2$ to $m=n$.  This yields that 
\begin{equation}
\label{Eq:MR23_01}
A_n \leq  |\delta_1| + |\delta_n| + \sum_{m=2}^{n} (c\tau_m + \tau_{m-1} + c\tau_m \tau_{m-1}) A_{m-1}	+  (\tau^{2(k+1)} + h^{2(r+1)}) (M_n + N_n) \,,
\end{equation}	
where by the definition of $\mathcal E_t^{n,i}$ and $\mathcal E_{\vec x}^{n,i}$, for $i\in \{1,2,3\}$,  there holds that 
\begin{subequations}
\label{Eq:MR24}	
\begin{alignat}{2}
	M_n :=\sum_{m=1}^n \big((\mathcal E_t^{n,1})^2 + (\mathcal E_t^{n,2})^2 + (\mathcal E_{\vec x}^{n,1})^2 + (\mathcal E_{\vec x}^{n,2})^2\big) \leq c< \infty\,, \quad N_n  :=\sum_{m=1}^n \tau \big((\mathcal E_{t}^{n,3})^2+(\mathcal E_{\vec x}^{n,3})^2 \big)\leq c < \infty
\end{alignat}
\end{subequations}	
for sufficiently regular solutions $(\vec u,p)$ to the system \eqref{Eq:HPS} and $n=1,\ldots,N$. We have that 
\begin{equation}
	\label{Eq:MR25}	
	\prod_{j=1}^{n-1} (1+c\tau_j) \leq \operatorname{e}^{cT} \,.
\end{equation}
Combining \eqref{Eq:DD_06}	and \eqref{Eq:MR23} yields that 
\begin{equation}
\label{Eq:MR25_01n1}	
|\delta_1| \leq c \tau_1^{2(k+1)}+ c h^{2(r+1)}\,.
\end{equation}
From the definitions \eqref{Eq:ES00}, \eqref{Eq:ErrP_00}, and \eqref{Eq:MR20} we conclude by the inequalities of Cauchy--Schwarz and Cauchy--Young and \eqref{Eq:Omega} that, for some sufficiently small $\varepsilon>0$,  there holds that 
\begin{equation}
\label{Eq:MR25_01}	
|\delta_n| = \langle \omega(t_n), \nabla \cdot \vec E_{\tau,h}^1(t_n) \rangle   \leq c h^{2(r+1)} + \varepsilon A_n
\end{equation}
The Gronwall argument, along with \eqref{Eq:MR24} to \eqref{Eq:MR25_01}	and Assumption \ref{AssIV}, then implies that
\begin{equation}
	\label{Eq:MR30}	
	||| \vec E_{\tau,h} (t_n)|||_e^2  + ||| e_{\tau,h}(t_n)|||_e^2  \leq c \tau^{2(k+1)}+ c h^{2(r+1)} \,, \quad \text{for}\;\; n= 0,\ldots, N\,,
\end{equation}
where $\tau = \max_{n=1,\ldots, N} \tau_n$: cf.~Subsec.~\ref{Subsec:FES}. 
By \eqref{Eq:ErrL2L201}, \eqref{Eq:En_1+01}, \eqref{Eq:MR30} and \eqref{Eq:Nequi} we then get that 
\begin{equation}
	\label{Eq:MR31}	
	||| \vec E_{\tau,h}|||_{L^2(I_n;\vec L^2)}^2  + \|e_{\tau,h}\|_{L^2(I_n;L^2)}^2  \leq c \tau \big(\tau^{2(k+1)}+c h^{2(r+1)}\big) 
\end{equation}
for $n=2,\ldots,N$. For $n=1$, estimate \eqref{Eq:MR31}	follows from \eqref{Eq:ErrL2L201} along with \eqref{Eq:MR12} and \eqref{Eq:Nequi}. By the $L^\infty$--$L^2$ inverse relation \eqref{Eq:LinfL2} we conclude from \eqref{Eq:MR31} that 
\begin{equation}
	\label{Eq:MR32}	
	||| \vec E_{\tau,h} (t)|||^2  + \|e_{\tau,h}(t)\|^2  \leq c \tau^{2(k+1)}+c h^{2(r+1)} \,, \quad \text{for}\;\; t\in [0,T] \,. 
\end{equation}
Finally, applying the triangle inequality to the splitting \eqref{Eq:ErrU_00} and \eqref{Eq:ErrP_00} of the errors  and employing the estimates \eqref{Eq:EtauEtav} proves the assertion \eqref{Eq:MR01}. For this, we note that \eqref{Eq:Etav} holds analogously for the error $\omega$ defined in \eqref{Eq:ErrP_00}; cf.\ \cite[Eq.\ (3.20)]{KM04}.	 
\end{mproof}

\begin{rem}
\label{Rem:OptOrder}
\begin{itemize}
\item We note that the constant of the error estimate \eqref{Eq:MR01} depends in particular on the norms of the continuous solution that are induced by Lem.~\ref{Lem:EstEtaOmega} and Lem.~\ref{Lem:ES} to Lem.~\ref{Lem:DifDelta}. Thereby, the tacitly assumed regularity of the continuous solution becomes obvious.

\item For arbitrray $t\in I$, estimate \eqref{Eq:MR01} is of optimal order with respect to the time and space discretization, if the approximation error is measured in terms of $\|p(t)-p_{\tau,h}(t)\|$ and the elastic energy quantity $\|\nabla (\vec u(t) -\vec u_{\tau,h}(t))\| + \| \vec v(t) -\vec v_{\tau,h}(t)\|$.

\item From \eqref{Eq:MR01}, an error estimate for $\|\vec u(t)-\vec u_{\tau,h}(t)\|$ can be obtained by the Poincar\'e inequality. However the resulting estimate for $\|\vec u(t)-\vec u_{\tau,h}\|$, as well as the estimate of $\|\vec v(t)-\vec v_{\tau,h}\|$ in \eqref{Eq:MR01}, are of suboptimal order with respect to the space discretization only. This is due to the coupling of the unknows of in the continuous system \eqref{Eq:HPS} and its fully discrete counterpart \eqref{Eq:DPQ_0}, the energy-type arguments of the error analysis bounding  the quantity $\|\nabla (\vec u(t) -\vec u_{\tau,h}(t))\| + \| \vec v(t) -\vec v_{\tau,h}(t)\|$ and, finally, the non-equal order approximation of $\vec u$ and $p$ by inf-sup stable pairs of finite element spaces. Similar observations regarding the coupling of the errors in the approximation of the unknowns are well-known from the discretization of the Navier--Stokes equations by inf-sup stable pairs of finite element spaces. In Sec.~\ref{Sec:NumExp}, the convergence rates of the error estimate \eqref{Eq:MR01} are confirmed by our numerical experiments. 

\item In \cite{KM04}, the convergence of a continuous Galerkin method for a scalar-valued nonlinear wave equation in $u$ is studied. Optimal order $L^2$-error estimates, for the quantities $u$ und $v=\partial_t u$, are proved. A key ingredient of this optimality is the special choice of the initial values, which is in contrast to our more general one given by Assumption~\ref{AssIV}. Compared to the purely hyperbolic case studied in  \cite{KM04}, in our analysis the projection error that is induced by the coupling term $\alpha \nabla \cdot \partial_t \vec u$ in \eqref{Eq:HPS_2} implies the loss of one order of accuracy for the spatial discretization of the overall system such that the result of \cite{KM04} regarding the $L^2$-error convergence of $u$ and $v$ cannot be transfered directly to the system \eqref{Eq:HPS}. Optimal order estimates for $\|\vec u - \vec u_{\tau,h}\|$ and $\|\vec v - \vec v_{\tau,h}\|$ might require proper decoupling techniques for the subproblems of \eqref{Eq:HPS} which has to be left as a work for the future.

\item We conjecture that the result \eqref{Intro:MainEst2} of superconvergence in the time nodes is satisfied. This is  illustrated numerically in Sec.~\ref{Sec:NumExp}. We expect that the proof of superconvergence can be built on Thm.\ \ref{Thm:MainRes}. However, this remains a work for the future. For the proof of superconvergence for the wave equation we refer to \cite{BKRS20}.

\end{itemize}
\end{rem}

\section{Numerical convergence test}
\label{Sec:NumExp}

Here we present the results of our performed numerical experiments in order to confirm Thm.~\ref{Thm:MainRes}.
The implementation of the numerical scheme was done in an in-house high-performance frontend solver for the \texttt{deal.II} library \cite{AB21}. We study \eqref{Eq:HPS} for $\Omega=(0,1)^2$ and $I=(1,2]$ and the prescribed solution
\begin{equation}
	\label{Eq:givensolution}
	\boldsymbol u(\boldsymbol x, t) = \phi(\boldsymbol x, t) \boldsymbol I_2 \;\; \text{and}\;\;
	p(\boldsymbol x, t) = \phi(\boldsymbol x, t) \;\; \text{with}\;\; 
	\phi(\boldsymbol x, t) = \sin(\omega_1 t^2) \sin(\omega_2 x_1) \sin(\omega_2 x_2)
\end{equation}
with $\omega_1=\omega_2 = \pi$. We put $\rho=1.0$, $\alpha=0.9$, $c_0=0.01$ and $\boldsymbol K=\boldsymbol I_2$ with the identity $\vec I_2\in \R^{2,2}$. For the fourth order elasticity tensor $\boldsymbol C$, isotropic material properties with Young's modulus $E=100$ and Poisson's ratio $\nu=0.35$ are chosen. In our experiments, the norm of $L^\infty(I;L^2)$ is approximated by ($t_{n,m}$: Gauss quadrature nodes of $I_n$)
\begin{equation*}
	\label{Eq:DiscNorm}
	\| w\|_{L^\infty(I;L^2)} \approx \max \{ \| w_{|I_n}(t_{n,m})\|  \mid m=1,\ldots ,M\,, \; n=1,\ldots,N\}\,, \quad \text{with}\;\; M=100\,.
\end{equation*}

We study the space-time convergence behavior of the scheme \eqref{Eq:DPQ_0} to confirm our main result \eqref{Eq:MR01}. For this, the domain $\Omega$ is decomposed into a sequence of successively refined meshes of quadrilateral finite elements. The spatial and temporal mesh sizes are halfened in each of the refinement steps. The step sizes of the coarsest space and time mesh are $h_0=1/(2\sqrt{2})$ and $\tau_0=0.05$. To illustrate  \eqref{Eq:MR01}, we choose the polynomial degree $k=2$ and $r=2$, such that discrete solutions $\vec u_{\tau,h}, \vec v_{\tau}\in (X_\tau^2(V_h^3))^2$ and $p_{\tau,h}\in X_\tau^2(V_h^2)$ are obtained, as well as $k=3$ and $r=3$ with $\vec u_{\tau,h}, \vec v_{\tau}\in (X_\tau^3(V_h^4))^2$ and $p_{\tau,h}\in X_\tau^3(V_h^3)$; cf.\ \eqref{Eq:DefXk} and \eqref{Eq:DefVh}. The calculated errors and corresponding experimental orders of convergence are summarized in Table~\ref{Tab:1} and \ref{Tab:2}, respectively. Table~\ref{Tab:1} and \ref{Tab:2} nicely confirm our main result \eqref{Eq:MR01}. The orders of convergence in time and space, expected from the estimate \eqref{Eq:MR01}, are clearly observed.

\begin{table}[H]
\centering
	\begin{tabular}{l}
	\begin{tabular}{cccccccc}
	\toprule
	{$\tau$} & {$h$} &
	{ $\| \nabla (\vec u - \vec u_{\tau,h})  \|_{L^2(\vec L^2)} $ } & {EOC} &
	{ $\| \vec v - \vec v_{\tau,h}  \|_{L^2(\vec L^2)} $ } & {EOC} &
	{ $\| p-p_{\tau,h}  \|_{L^2(L^2)}  $ } & {EOC} 
 \\
	\cmidrule(r){1-2}
	\cmidrule(lr){3-8}
	$\tau_0/2^0$ & $h_0/2^0$ & 3.7772346728e-03 & {--} & 4.4831153608e-03  & {--} & 1.3925593715e-03  & {--} \\ 
	$\tau_0/2^1$ & $h_0/2^1$ & 4.7293499671e-04 & 3.00 & 5.6200459009e-04  & 3.00 & 1.7624666295e-04  & 2.98 \\
	$\tau_0/2^2$ & $h_0/2^2$ & 5.9118396929e-05 & 3.00 & 7.0409147572e-05  & 3.00 & 2.2094955372e-05  & 3.00\\
	$\tau_0/2^3$ & $h_0/2^3$ & 7.3894810579e-06 & 3.00 & 8.8070050157e-06  & 3.00 & 2.7638964740e-06  & 3.00\\
	\bottomrule
\end{tabular}\\
\mbox{}\\
	\begin{tabular}{cccccccc}
	\toprule
	{$\tau$} & {$h$} &
	{ $\| \nabla (\vec u - \vec u_{\tau,h})   \|_{L^{\infty}(\vec L^2)} $ } & {EOC} &
	{ $\| \vec v - \vec v_{\tau,h}  \|_{L^{\infty}(\vec L^2)} $ } & {EOC} &
	{ $\| p-p_{\tau,h}  \|_{L^{\infty}(L^2)} $ } & {EOC} \\
	\cmidrule(r){1-2}
	\cmidrule(lr){3-8}
	$\tau_0/2^0$ & $h_0/2^0$ & 5.5609986126e-03 & {--} & 1.4388258226e-02 & {--} & 1.9457909519e-03 & {--}  \\ 
	$\tau_0/2^1$ & $h_0/2^1$ & 7.3872532490e-04 & 2.91 & 1.8026863849e-03 & 3.00 & 2.4740005168e-04 & 2.98  \\
	$\tau_0/2^2$ & $h_0/2^2$ & 9.4556857326e-05 & 2.97 & 2.2667403592e-04 & 2.99 & 3.0867702485e-05 & 3.00  \\
	$\tau_0/2^3$ & $h_0/2^3$ & 1.1925250119e-05 & 2.99 & 2.8448677188e-05 & 2.99 & 3.8601048383e-06 & 3.00  \\
	\bottomrule
\end{tabular}
\end{tabular}
\caption{%
		$L^2(L^2)$ and $L^\infty(L^2)$ errors and experimental orders of convergence (EOC) for \eqref{Eq:givensolution} with polynomial degrees $k=2$ and $r=2$.
	}
	\label{Tab:1}
\end{table}

\begin{table}[H]
	\centering
	\begin{tabular}{l}
	\begin{tabular}{cccccccc}
				\toprule
	{$\tau$} & {$h$} &
	{ $\| \nabla (\vec u - \vec u_{\tau,h})  \|_{L^2(\vec L^2)} $ } & {EOC} &
	{ $\| \vec v - \vec v_{\tau,h}\|_{L^2(\vec L^2)}  $ } & {EOC} &
	{ $\| p - p_{\tau,h} \|_{L^2(L^2)}  $ } & {EOC}  \\
	\cmidrule(r){1-2}
	\cmidrule(lr){3-8}
	$\tau_0/2^0$ & $h_0/2^0$ & 1.7724800037e-04 & {--} & 1.5572598126e-04  & {--} & 6.2865996817e-05  & {--} \\ 
	$\tau_0/2^1$ & $h_0/2^1$ & 1.1068826736e-05 & 4.00 & 9.0324299079e-06  & 4.11 & 3.9664381213e-06  & 3.99  \\
	$\tau_0/2^2$ & $h_0/2^2$ & 6.9153355647e-07 & 4.00 & 5.5554036618e-07  & 4.02 & 2.4851816029e-07  & 4.00\\
	$\tau_0/2^3$ & $h_0/2^3$ & 4.3215752542e-08 & 4.00 & 3.4586146527e-08  & 4.01 & 1.5542077250e-08  & 4.00  \\
	\bottomrule
	\end{tabular}\\
\mbox{}\\
\begin{tabular}{cccccccc}
\toprule
{$\tau$} & {$h$} &
	{ $\| \nabla (\vec u - \vec u_{\tau,h})  \|_{L^\infty(\vec L^2)} $ } & {EOC} &
{ $\| \vec v - \vec v_{\tau,h}\|_{L^\infty(\vec L^2)}  $ } & {EOC} &
{ $\| p - p_{\tau,h} \|_{L^\infty(L^2)}  $ } & {EOC}  \\
\cmidrule(r){1-2}
\cmidrule(lr){3-8}
$\tau_0/2^0$ & $h_0/2^0$ &  3.0383309559e-04 & {--} & 5.7065321892e-04 & {--} & 9.3580580659e-05 & {--}  \\ 
$\tau_0/2^1$ & $h_0/2^1$ &   1.9175723302e-05 & 3.99 & 3.8885259584e-05 & 3.88 & 5.8271904381e-06 & 4.01  \\
$\tau_0/2^2$ & $h_0/2^2$ &  1.1977037979e-06 & 4.00 & 2.5396723780e-06 & 3.94 & 3.6728075814e-07 & 3.99  \\
$\tau_0/2^3$ & $h_0/2^3$ &  7.4962458146e-08 & 4.00 & 1.6227333767e-07 & 3.97 & 2.3002686673e-08 & 4.00  \\
\bottomrule	
\end{tabular}	
     \end{tabular}

	\caption{%
		$L^2(L^2)$ and $L^\infty(L^2)$  errors and experimental orders of convergence (EOC) for \eqref{Eq:givensolution} with polynomial degrees $k=3$ and $r=3$.	
	}
	\label{Tab:2}
\end{table}

 In Table~\ref{Tab:3}, superconvergence in the discrete time nodes is studied in terms of the time mesh dependent norm
 \begin{equation}
 \label{Eq:MDN}
 	\| w\|_{l^\infty(L^2)} := \max\{ \| w(t_n)\| \mid n=1,\ldots, N\}\,.
 \end{equation}
  For the finite element spaces we choose the orders $k=3$ and $r=5$ such that discrete solutions $\vec u_{\tau,h}, \vec v_{\tau}\in (X_\tau^3(V_h^6))^2$ and $p_{\tau,h}\in X_\tau^3(V_h^5)$ are obtained. Superconvergence of order $2k$ in the discrete time nodes is clearly observed in Table~\ref{Tab:3}. This confirms our conjecture~\eqref{Intro:MainEst2} of superconvergence of order $2k$ in the discrete time nodes $t_n$, for $n=1,\ldots,N$. 
\begin{table}[H]
\centering
\begin{tabular}{l}
\begin{tabular}{cccccccc}
\toprule
{$\tau$} & {$h$} &
{ $\| \nabla (\vec u - \vec u_{\tau,h})  \|_{L^2(\vec L^2)} $ } & {EOC} &
{ $\| \vec v - \vec v_{\tau,h}\|_{L^2(\vec L^2)}  $ } & {EOC} &
{ $\| p - p_{\tau,h} \|_{L^2(L^2)}  $ } & {EOC}  \\
\cmidrule(r){1-2}
\cmidrule(lr){3-8}
$\tau_0/2^0$ & $h_0/2^0$ & 5.8117734426e-05 & {--} & 1.5347090551e-04  & {--} & 9.3413974336e-06  & {--} \\ 
$\tau_0/2^1$ & $h_0/2^1$ & 3.6198825671e-06 & 4.00 & 8.9954777890e-06  & 4.09 & 5.7613608543e-07  & 4.02 \\
$\tau_0/2^2$ & $h_0/2^2$ & 2.2603227629e-07 & 4.00 & 5.5496215896e-07  & 4.02 & 3.5977539073e-08  & 4.00   \\
$\tau_0/2^3$ & $h_0/2^3$ & 1.4123671689e-08 & 4.00 & 3.4577094422e-08  & 4.00 & 2.2483070160e-09  & 4.00   \\
\bottomrule	
\end{tabular}\\
\mbox{}\\

\begin{tabular}{cccccccc}
\toprule
{$\tau$} & {$h$} &
{ $\| \nabla (\vec u - \vec u_{\tau,h})  \|_{l^\infty(\vec L^2)} $ } & {EOC} &
{ $\| \vec v - \vec v_{\tau,h}\|_{l^\infty(\vec L^2)}  $ } & {EOC} &
{ $\| p - p_{\tau,h} \|_{l^\infty(L^2)}  $ } & {EOC}  \\
\cmidrule(r){1-2}
\cmidrule(lr){3-8}
$\tau_0/2^0$ & $h_0/2^0$ & 1.1089049623e-05 & {--} & 1.4804895672e-04 & {--} & 1.0389805110e-05 & {--}  \\ 
$\tau_0/2^1$ & $h_0/2^1$ &  1.4735513623e-07 & 6.23 & 2.1095147908e-06 & 6.13 & 1.2944103974e-07 & 6.33  \\
$\tau_0/2^2$ & $h_0/2^2$ &  2.3655340792e-09 & 5.96 & 3.3680209502e-08 & 5.97 & 2.1560790646e-09 & 5.91  \\
$\tau_0/2^3$ & $h_0/2^3$ & 3.6038421330e-11 & 6.04 & 5.2092447939e-10 & 6.01 & 3.3031484852e-11 & 6.03  \\
\bottomrule	
\end{tabular}
\end{tabular}

	\caption{$L^2(L^2)$ and $l^\infty(L^2)$  errors (cf.\ \eqref{Eq:MDN}) and experimental orders of convergence (EOC) for \eqref{Eq:givensolution} with polynomial degrees $k=3$ and $r=5$, showing superconvergence in the discrete time nodes. 	
	}
	\label{Tab:3}
\end{table}

\begin{rem}[Iterative solver for the algebraic system]
Higher order variational time discretizations, corresponding to larger values of the polynomial degree $k$, lead to complex block matrices on the algebraic level. Their efficient iterative solution is a challenging task. For this, we use GMRES iterations that are preconditioned by a $V$-cycle of the geometric multigrid method. For the smoothing operations a local Vanka method is applied. For further details of the design of the solver for space-time finite element methods and numerical experiments demonstrating its efficiency and robustness we refer to our work \cite{AB22,AB21} on the application of such techniques to the Navier--Stokes system. The presentation and numerical study of the geometric multigrid preconditioner for the dynamic Biot system \eqref{Eq:HPS} as well as three-dimensional simulations of the Biot system will be addressed in a forthcoming work.   
\end{rem}

\section*{Acknowledgement}

This work was supported by the German Academic Exchange Service (DAAD) under the grant ID 57458510 and by the Research Council of Norway (RCN) under the grant ID 294716.  F.\ A.\ Radu acknowledges funding from the VISTA programme, The Norwegian Academy of Science and Letters.


\begin{thebibliography}{99}
	
\bibitem{AMN09}
G.\ Akrivis, C.\ Makridakis, R.\ H.\ Nochetto, {\em Optimal order a posteriori 
	error estimates for a class of Runge--Kutta and Galerkin methods}, Numer.\ 
Math., \textbf{114} (2009), pp.\ 133--160.

\bibitem{AMN11}
G.\ Akrivis, C.\ Makridakis, R.\ H.\ Nochetto, {\em Galerkin and Runge--Kutta 
	methods: unified formulation, a posteriori error estimates and nodal
	superconvergence}, Numer.\ Math., \textbf{118} (2011), pp.\ 429--456.

\bibitem{AB22}
M.\ Anselmann, M.\ Bause, \textit{Efficiency of local Vanka smoother geometric multigrid preconditioning for space-time finite element methods to the Navier–Stokes equations}, PAMM Proc.\ Appl.\ Math.\ Mech., \textbf{accepted} (2022), pp.\ 1--6; arXiv:2210.02690. 

\bibitem{AB21}
M.\ Anselmann, M.\ Bause, \textit{A geometric multigrid method for space-time finite element
discretizations of the Navier–Stokes equations and its application to 3d flow simulation},
ACM Trans.\ Math.\ Softw., \textbf{submitted} (2021), pp. 1--27; arXiv:2107.10561.

\bibitem{AB20}
M. Anselmann, M. Bause,  Higher order Galerkin-collocation time discretization with Nitsche’s method for the Navier-Stokes equations,  Math. Comp. Simul., in press (2020), pp.\ 1--22; doi:10.1016/j.matcom.2020.10.027.

\bibitem{ABBM20}
M.\ Anselmann, M.\ Bause, S.\ Becher, G.\ Matthies, \textit{Galerkin–collocation approximation in time for the wave equation and its post-processing}, ESAIM: M2AN, \textbf{54} (2020), pp.~2099--2123.  

\bibitem{ADMQ16}
P.\ F.\ Antonietti, B.\ A.\ De Dios, I.\ Mazzieri, A.\ Quarteroni, \textit{Stability analysis of discontinuous Galerkin approximations to the elastodynamics problem}, J.\ Sci.\ Comput., \textbf{68} (2016), pp.~143–170.

\bibitem{ABB21}
D.\ Arndt, W.\ Bangerth, B.\ Blais, M.\ Fehling, R.\ Gassmöller, T.\ Heister, L.\ Heltai, U.\ Köcher, M.\ Kronbichler, M.\ Maier, P.\ Munch, J.\ Pelteret, S.\ Proell, K.\ Simon, B.\ Turcksin, D.\ Wells, J.\ Zhang, J. \textit{The deal.II Library, Version 9.3}, J.\ Numer.\ Math., .\textbf{29} (2021), pp.\ 171--186.

\bibitem{AM89}
A.\ K.\ Aziz, P.\ Monk, \textit{Continuous finite elements in space and time for the heat equation}, Math.\ Comp., \textbf{52} (1989), pp.\ 255--274. 


\bibitem{B19}
M.\ Bause, \textit{Iterative coupling of mixed and discontinuous Galerkin  methods for poroelasticity}, in F.\ A.\ Radu et al. (Eds.), \emph{Numerical Mathematics and Advanced Applications ENUMATH 2017}, Lecture Notes in Computational Science and Engineering 126, Springer, Cham, 2019, pp.\ 551--560.

\bibitem{BKRS20}
M.\ Bause, U.\ K\"ocher, F.\ A.\ Radu, F.\ Schieweck, \textit{Post-processed Galerkin approximation of improved order for wave equations}, Math.\ Comp., \textbf{89} (2020), pp.\ 595--627.

\bibitem{BRK17}
M.\ Bause, F.\ A.\ Radu, U.\ K\"ocher,  \textit{Space-time finite element approximation of the Biot poroelasticity system with iterative coupling}, Comput.\ Methods Appl.\ Mech.\ Engrg., \textbf{320} (2017), pp.\ 745--768.

\bibitem{BM21}
S.\ Becher and G.\ Matthies, \textit{Variational time discretizations of higher order and higher regularity}, BIT Numer.\ Math., \textbf{61} (2021), pp.\ 721--755.

\bibitem{B41}
M.\ Biot, \textit{General theory of three-dimensional consolidation}, J.\ Appl.\ Phys., \textbf{12} (1941), pp.~155--164.

\bibitem{B55}
M.\ Biot, \textit{Theory of elasticity and consolidation for a porous anisotropic solid}, J.\ Appl.\ Phys., \textbf{26} (1955), pp.~182--185.

\bibitem{B72}
M.\ Biot, \textit{Theory of finite deformations of porous solids, Indiana Univ.\ Math.\ J.}, \textbf{21} (1972), pp.~597--620.

\bibitem{BKNR22}
J.\ W.\ Both, N.\ A.\ Barnafi, F.\ A.\ Radu, P.\ Zunino, A.\ Quarteroni, \textit{Iterative splitting schemes for a soft material poromechanics model}, Comput.\ Methods Appl.\ Mech.\ Engrg., \textbf{388} (2022), 114183.

\bibitem{BS94}
S.\ C.\ Brenner, L.\ R.\ Scott, \textit{The Mathematical Theory of Finite Element Methods}, Springer, New York, 1994. 

\bibitem{C72}
D.\ E.\ Carlson, \textit{Linear thermoelasticity}, Handbuch der Physik  V Ia/2, Springer, Berlin, 1972.

\bibitem{CW81}
S.-I.\ Chou,  C.-C.\ Wang, \textit{ Estimates of error in finite element approximate solutions to problems in linear thermoelasticity. Part 1. Computationally coupled numerical schemes}, Arch.\ Rational Mech.\ Analy, \textbf{77} (1981), pp.~263--299.

\bibitem{PE12}
D.\ A.\ Di Pietro, A.\ Ern,  \textit{Mathematical Aspects of Discontinuous Galerkin Methods}, Springer, Heidelberg, 2012.  

\bibitem{DF15}
V.\  Dolej\v{s}\'{i}, M.\ Feistauer, \textit{Discontinuous Galerkin Method}, Springer, Heidelberg, 2015. 

\bibitem{FP96}
D.\ A.\ French, T.\ E.\ Peterson,  \textit{A continuous space-time finite element method for the wave equation}, Math.\ Comp., \textbf{65} (1996), pp.\ 491--506.

\bibitem{F47}
K.\ O.\ Friedrichs, \textit{On the boundary value problems of the theory of elasticity and Korn's inequality}, Ann.\ of Math., \textbf{48} (1947), pp.\ 441--471.

\bibitem{GSS06} 
M.\ J.\ Grote, A.\ Schneebeli, D.\ Sch{\"{o}}tzau, \textit{Discontinuous {G}alerkin finite element method for the wave equation}, SIAM J.\ Numer.\ Anal., \textbf{44} (2006), pp.\ 2408--2431.

\bibitem{GS09} 
M.\ J.\ Grote, D.\ Sch\"otzau, \textit{Optimal error estimates for
the fully discrete interior penalty {DG} method for the wave equation},
J.\ Sci.\ Comput., \textbf{40} (2009), pp.\ 257--272.

\bibitem{HST13}
S.\ Hussain, F.\ Schieweck, S.\ Turek, \textit{An efficient and stable finite element solver of higher order in space and time for nonstationary incompressible flow}, Internat.\ J. Numer.\ Methods Fluids, \textbf{73} (2013), pp.\ 927--952.

\bibitem{H91}
G.\ W.\ Howell, \textit{Derivative bound for Lagrange interpolation: An extension of Cauchy's bound for error of Lagrange interpolation}, J. Approx.\ Theory, \textbf{67} (1991), pp.\ 164--173.

\bibitem{HH88}
T.\ J.\ R.\ Hughes, G.\ M.\ Hulbert, {\em Space-time finite element methods for  elastodynamics: Formulations and error estimates}, Comput.\ Methods Appl.\ Mech.\ Engrg., 
\textbf{66} (1988), pp.\ 339--363.

\bibitem{H72_1}
B.\ L.\ Hulme, \emph{One-step piecewise polynomial Galerkin methods for initial value problems}, Math.\ Comp., \textbf{26} (1972), pp.\ 416--426.

\bibitem{H72_2}
B.\ L.\ Hulme, \emph{Discrete Galerkin and related one-step methods for ordinary differential equations}, Math.\ Comp., \textbf{26} (1972), pp.\ 881--891.

\bibitem{JR18}
S.\ Jiang, R.\ Racke, \textit{Evolution equations in thermoelasticity}, CRC Press, Boca Raton, 2018.

\bibitem{J16}
V.\ John, \emph{Finite Element methods for Incompressible Flow Problems}, Springer, Cham 
2016.

\bibitem{J93}
C.\ Johnson, {\em Discontinuous Galerkin finite element methods for second order hyperbolic problems}, Comput.\ Methods Appl.\ Mech.\ Engrg., \textbf{107} (1993), pp.\ 117--129.

\bibitem{KM99}
O.\ Karakashian, C.\ Makridakis, \textit{A space-time finite element method for the nonlinear Schrödinger equation: The continuous Galerkin method}, SIAM J.\ Numer.\ Anal., \textbf{36} (1999), pp.\ 1779--1807.

\bibitem{KM04}
O.\ Karakashian, C.\ Makridakis, \textit{Convergence of a continuous Galerkin method with mesh modification for nonlinear wave equations}, Math.\ Comp., \textbf{74} (2004), pp.\ 85--102.

\bibitem{K15}
U.\ K\"ocher, \textit{Variational space-time methods for the elastic wave equation and the diffusion equation}, PhD Thesis, Helmut-Schmidt-Universit\"at,  \url{http://edoc.sub.uni-hamburg.de/hsu/volltexte/2015/3112/}, 2015.
\bibitem{KB14}

U.\ K\"ocher, M.\ Bause, \textit{Variational space-time methods for the wave equation}, J.\ 
Sci.\ Comput., \textbf{61} (2014), pp.\ 424--453.

\bibitem{LLW16}
S.\ Lee, Y.\ J.\ Lee, M.\ F.\ Wheeler, \textit{A locally conservative enriched Galerkin approximation and efficient solver for elliptic and parabolic problems}, SIAM J.\ Sci.\ Comput., \textbf{38} (2016), pp.\ A1404-A1429.

\bibitem{L86}
R.\ Leis, \textit{Initial boundary value problems in mathematical physics}, Teubner, Stuttgart, John Wiley \& Sons, Chichester, 1986.

\bibitem{MW12}
A.\ Mikeli\'{c}, M.\ F.\ Wheeler, {\em Theory of the dynamic Biot--Allard equations and their link to the quasi-static Biot system}, J.\ Math.\ Phys., \textbf{53} (2012), 123702:1--15.

\bibitem{MW13}
A.\ Mikeli\'{c}, M.\ F.\ Wheeler, {\em Convergence of iterative coupling for coupled flow and geomechanics}, Comput.\  Geosci., \textbf{17} (2013), pp.\ 479--496.

\bibitem{ML92}
M.\ A.\ Murad, A.\ F.\ D.\ Loula, \textit{Improved accuracy in finite element analysis of Biot’s consolidation problem}, Comput.\ Methods Appl.\ Mech. Engrg., \textbf{95} (1992), pp.~359--382.

\bibitem{ML94}
M.\ A.\ Murad, A.\ F.\ D.\ Loula, \textit{On stability and convergence of finite element approximations of Biot’s consolidation problem}, Internat. J.\ Numer.\ Methods Engrg., \textbf{37} (1994), pp.~645--667.

\bibitem{MTL96}
M.\ A.\ Murad, V.\ Thom\'{e}e, A.\ F.\ D.\ Loula, \textit{Asymptotic behavior of semidiscrete finite-element approximations of Biot’s consolidation problem}, SIAM J.\ Numer.\ Anal., \textbf{33} (1996), pp.~1065--1083.

\bibitem{ORSZ21}
R.\ Oyarz\'{u}a, S.\ Rhebergen, M.\ Solano, P.\ Z\'{u}\~{n}iga, \textit{Error analysis of a conforming and locking-free four-field formulation for the stationary Biot's model}, ESAIM: M2AN, \textbf{55} (2021), pp.~475--506.

\bibitem{PW08}
P.\ Philips, M.\ Wheeler, \textit{A coupling of mixed and discontinuous Galerkin finite element methods for poroelasticity}, Comput.\ Geosci., \textbf{12} (2008),  pp.~417--435.

\bibitem{QV94}
A.\ Quarteroni, V.\ Valli, \textit{Numerical Approximation of Partial Differential Equations}, Springer, Berlin, 1994.  

\bibitem{RHOAGZ17}
C.\ Rodrigo, X.\ Hu, P.\ Ohm, J.\ H.\ Adler, F.\ J.\ Gaspar, L.\ T.\ Zikatanov, \textit{New stabilized discretizations for poroelasticity and the Stokes' equations}, Comput.\ Methods Appl.\ Mech.\ Engrg., \textbf{341} (2018), pp.\ 467--484.

\bibitem{S10}
F.\ Schieweck, \textit{A-stable discontinuous Galerkin–Petrov time discretization of higher order}, J.\ Numer.\ Math., \textbf{18} (2010), pp.\ 25--57.

\bibitem{STW22}
C.\ Seifert, S.\ Trostorff, M.\ Waurick, \textit{Evolutionary Equations: Picard's Theorem for Partial Differential Equations, and Applications}, Birkhäuser, Cham, 2022.

\bibitem{S00}
R.\ Showalter, \textit{Diffusion in poro-elastic media}, J.\ Math.\ Anal.\ Appl.,  \textbf{251} (2000), pp.\ 310--340.

\bibitem{SZ20}
O.\ Steinbach, M.\ Zank, \textit{Coercive space-time finite element methods for initial boundary value problems}, Electron.\ Trans.\ Numer.\ Anal., \textbf{52} (2020), pp.\ 154--194.

\bibitem{S89}
M.\ Slodi\v{c}ka, \textit{Application of Rothe's method to integrodifferential equation}, Comment.\ Math.\ Univ.\ Carolinae, \textbf{30} (1989), pp.~57--70.

\bibitem{SZ22}
O.\ Steinbach, M.\ Zank, \textit{A generalized inf–sup stable variational formulation for the wave equation}, J.\ Math.\ Anal.\ Appl.,  \textbf{505} (2022), 125457. 

\bibitem{SL09}
S.\ Sun, J.\ Liu, \textit{A locally conservative finite element method based on piecewise constant enrichment of the continuous Galerkin method}, SIAM J.\ Sci.\ Comput., \textbf{31} (2009), pp.\ 2528--2548.

\bibitem{T06}
V. Thome{\'{e}}, {\em Galerkin Finite Element Methods for Parabolic Problems}, Springer, Berlin, 2006.

\bibitem{VSBW18}
J.\ Vamaraju, M.\ K.\ Sen, J.\ De Basabe, M.\ Wheeler,
\textit{Enriched Galerkin finite element approximation for elastic wave propagation in fractured media}, J.\ Comput.\ Phys., \textbf{372} (2018), pp.\ 726--747.

\bibitem{VR15}
M.\ Vlasak, F.\ Roskovec, \emph{On Runge-Kutta, collocation and discontinuous Galerkin methods: Mutual connections and resulting consequences to the analysis}, Programs and Algorithms of Numerical Mathematics, \textbf{17} (2015), pp.\ 231--236.

\bibitem{W16}
J.\ A.\ White, N.\ Castelletto, H.\ A.\ Tchelepi, {\em Block-partitioned solvers for coupled poromechanics: A unified framework}, Comp.\ Meth.\ Appl.\ Mech.\ Eng., 
\textbf{303} (2016), pp.\ 55--74.

\end{thebibliography}
\end{document}